\documentclass[12pt]{article}
\usepackage{amssymb}
\usepackage{tikz}
\usepackage{mathrsfs}
\font\eufm=eufm10 at 14pt\font\eufms=eufm10\font\eufmss=eufm7\newfam\eufam
\textfont\eufam=\eufm\scriptfont\eufam=\eufms\scriptscriptfont\eufam=\eufmss
\def\build#1_#2^#3{\mathrel{\mathop{\kern 0pt#1}\limits_{#2}^{#3}}}
\def\Z{{\bf Z}}\def\C{{\bf C}}\def\H{{\rm H}}\def\Hom{{\rm Hom}}\def\fl{\longrightarrow}
\def\Map{\hbox{Map}}\def\colim{\hbox{colim\ }}
\def\rde{\mathscr}
\def\C{{\rde C}}

\def\cqfd{\hfill\vbox{\hrule\hbox{\vrule height6pt depth0pt\hskip 6pt \vrule height6pt}\hrule\relax}}
\def\noi{\noindent}
\def\hfl#1#2{\smash{\mathop{\hbox to 12 mm{\rightarrowfill}}\limits^{\scriptstyle#1}_{\scriptstyle#2}}}
\def\vfl#1#2{\llap{$\scriptstyle #1$}\left\downarrow\vbox to 6mm{}\right.\rlap{$\scriptstyle #2$}}
\def\mfl#1#2{\smash{\mathop{\hbox to 10 mm{\leftarrowfill}}\limits^{\scriptstyle#1}_{\scriptstyle#2}}}
\def\diagram#1{\def\normalbaselines{\baselineskip=0pt\lineskip=10pt\lineskiplimit=1pt} \matrix{#1}}
\def\pv{\raise 2pt\hbox{$\bigwedge$}}
\begin{document}

\overfullrule=0pt
\
\vskip 64pt
\centerline{\bf Codimension $2$ embeddings and finite localization of spaces.}
\vskip 12pt
\centerline{Pierre Vogel\footnote{Universit\'e Paris Diderot, Institut de Math\'ematiques de Jussieu-Paris Rive Gauche (UMR 7586),
B\^atiment Sophie Germain, Case 7012, 75205--Paris Cedex 13 --- Email: pierre.vogel@imj-prg.fr}}
\vskip 48pt
\noi{\bf Abstract.} In order to classify concordance classes of codimension $2$ embeddings in a manifold $M$, we need to determine the complement
of such an embedding. These complements are spaces over $M$ well defined up to some homology equivalence. We construct a localization functor corresponding to this class
of homology equivalences and we give a characterization of local objects in terms of homotopy groups.

\vskip 12pt
\noi{\bf Keywords:} Non-commutative localization, localization of spaces.

\noi{\bf Mathematics Subject Classification (2020):} 18E35, 55P60
\vskip 24pt
\noi{\bf Introduction.}
A very important problem in topology is the problem of codimension $2$ embeddings. If we want to classify the concordance classes of embeddings from a closed manifold $V$ 
to a connected closed manifold $M$, we need to determine the complement $X$ of this embedding. The homology of $X$ with coefficients in $\Z[\pi_1(M)]$ can be determined 
by the geometry of the embedding. If the codimension of the embedding is more than $2$, the fundamental group of $X$ is isomorphic to $\pi_1(M)$ and the homotopy type of 
$X$ can be determined (in some sense). But, in the codimension $2$ case, the situation is completely different because the morphism $\pi_1(X)\rightarrow\pi_1(M)$ is never
an isomorphism.

If two embeddings $i_0:V\fl M$ and $i_1:V\fl M$ are concordant, we have the complements $X_0$ and $X_1$ of these embeddings and the complement $Y$ of the concordance
$V\times[0,1]\fl M\times[0,1]$. Up to homotopy, the only thing we know is the fact that the maps $X_0\subset Y$ and $X_1\subset Y$ induce isomorphisms in homology with
coefficients in $\Z[\pi_1(M)]$.

The goal of this paper is to give a good description of some localization functor obtained by inverting, in a suitable category, these kind of maps.

Consider a topological space $B$. By a $B$-space (or a space over $B$) we mean a pair $(X,f)$ where $X$ is a space having the homotopy type of a CW complex and $f:X\fl B$
is a continuous map. A $B$-map $\varphi:(X,f)\fl (Y,g)$ is a continuous map $\varphi:X\fl Y$ such that $g\varphi=f$. So we get a category: the category of $B$-spaces
denoted by Sp$_B$.

If $B$ is pointed, we say that a $B$-space $(X,f)$ is pointed if $X$ is pointed and $f$ preserves basepoints.

We say that a morphism $\varphi:(X,f)\fl (Y,g)$ is a cofibration or a homotopy equivalence if the induced map $X\fl Y$ is a cofibration or a homotopy equivalence.

In the previous situation, the complement of an embedding or a concordance in a manifold $M$ may be considered as a $M$-space.

Suppose $B$ is a pointed path-connected space with fundamental group $\pi$ and $\Z[\pi]\fl R$ is a ring homomorphism. In this situation, we have a class $W(B,R)$ of
morphisms in Sp$_B$ defined as follows:

A morphism $\varphi:(X,f)\fl (Y,g)$ belongs to $W(B,R)$ if and only if $X$ and $Y$ are finite CW complexes and $\varphi:X\fl Y$ is an inclusion inducing a bijection
$\pi_0(X)\build\fl_{}^\sim\pi_0(Y)$ and an isomorphism $\H_*(X,R)\build\fl_{}^\sim\H_*(Y,R)$.

If $X_0\fl Y$ and $X_1\fl Y$ are induced by a concordance of embeddings in $M$, these maps are in the class $W(M,\Z[\pi])$ (up to homotopy).
\vskip 12pt
Consider a category $\C$ and a class of morphisms in $\C$ denoted by $W$. An object $X$ in $\C$ is said to be $W$-local if, for every morphism $f:Y\fl Z$ in $W$, the 
induced map $f^*:\Hom(Z,X)\fl\Hom(Y,X)$ is bijective and a morphism $f:Y\fl Z$ is said to be a $W$-equivalence if, for every $W$-local object $X$, the map 
$f^*:\Hom(Z,X)\fl\Hom(Y,X)$ is bijective. A morphism $f:Y\fl Z$ is said to be a $W$-localization if $f$ is a $W$-equivalence and $Z$ is $W$-local.

A localization in $\C$ with respect to $W$ is a pair $(\Phi,i)$ where $\Phi:\C\fl\C$ is a functor and $i:$Id$\fl\Phi$ is a morphism of functors such that, for every 
$X$ in $\C$, $i:X\fl \Phi(X)$ is a $W$-localization. In this situation, we can see that an object $X\in\C$ is $W$-local if and only if $i:X\fl\Phi(X)$ is an isomorphism.
\vskip 12pt
For every $B$-space $(X,f)$ we may set: $(X,f)\times[0,1]=(X\times[0,1],f\circ pr_1)$ and we get a well defined notion of homotopy in Sp$_B$. Using that, we say that a
map $\varphi:U\fl V$ in Sp$_B$ is a strong homotopy equivalence if there is a map $\psi:V\fl U$ such that $\varphi\psi$ and $\psi\varphi$ are homotopic to identities.

In the category Sp$_B$ we want to inverse morphisms in a class $W$ only up to homotopy and we have to modify the notions of $W$-local and localization in the following
way:

Consider two $B$-spaces $X$ and $Y$. With the compact-open topology, the set of maps $X\fl Y$ is a topological space denoted by $\Map(X,Y)$. 
Suppose $W$ is a class of cofibrations in Sp$_B$. We say that a $B$-space $X$ is $W$-local if, for every $f: Y\fl Z$ in $W$, the induced map $f^*:\Map(Z,X)\fl\Map(Y,X)$
is a Serre fibration with weakly contractible fibers (i.e. such that each fiber has the weak homotopy type of a point). We say also that a map $f:Y\fl Z$ in Sp$_B$ is a
$W$-equivalence if $f$ is a cofibration and, for every $W$-local $B$-space $X$, the induced map $f^*:\Map(Z,X)\fl\Map(Y,X)$ is a Serre fibration with weakly contractible 
fibers. As above a morphism $f:Y\fl Z$ in Sp$_B$ is said to be a $W$-localization if $f$ is a $W$-equivalence and $Z$ is $W$-local.

A localization in Sp$_B$ with respect to $W$ is a pair $(\Phi,i)$ where $\Phi:$Sp$_B\fl$ Sp$_B$ is a functor and $i:$Id$\fl\Phi$ is a morphism of functors such that for 
every $B$-space $X$, $i:X\fl \Phi(X)$ is a $W$-localization. As above, we can check that a $B$-space $X$ is $W$-local if and only if
$i:X\fl \Phi(X)$ is a strong deformation retract. 

For any category $\C$ we let Mor$(\C)$ denote the category whose objects are the morphisms of $\C$ and whose morphisms are the commutative squares in $\C$.

Let us denote by $\overline W$ the class of cofibration $\varphi$ in Sp$_B$ such that, for every inclusion $(K\subset L)$
of finite $B$-complexes, every morphism $(K\subset L)\fl\varphi$ in Mor$($Sp$_B)$ factors through a morphism in $W$.

\vskip 12pt
\noi{\bf Theorem 1:} {\sl Let $B$ be a pointed path-connected space with fundamental group $\pi$, $R$ be a ring and $\Z[\pi]\fl R$ be a ring homomorphism. Then there
  exists a localization $(L^B_R,i)$ in Sp$_B$ with respect to the class $W(B,R)$.

Moreover, for every $B$-space $X$, the morphism $i:X\fl L^B_R(X)$ belongs to the class $\overline W(B,R)$. }
\vskip 12pt
Since $L^B_R$ is the localization with respect to $R$-homology equivalences between finite complexes, $(L^B_R,i)$ will be called a finite localization in Sp$_B$ with
respect to $R$-homology.
\vskip 12pt
Using this theorem, we get a strategy to study the concordance classes of embeddings $V\subset M$. We may replace a complement of an embedding by its localization and we
get a homotopy cocartesian diagram:
$$\diagram{S(\xi)&\hfl{}{}&\widehat X\cr\vfl{}{}&&\vfl{}{}\cr B(\xi)&\hfl{}{}&M\cr}$$
where $\xi$ is a vector bundle on $V$, $S(\xi)$ and $B(\xi)$ are total spaces of their corresponding sphere and disk bundles and $\widehat X$ is a
$W(M,\Z[\pi_1(M)])$-local $M$-space. Moreover this diagram depends only on the concordance class of the given embedding.

Conversely, if we have such a homotopy cocartesian diagram, we have to replace $\widehat X$ by a manifold and reconstruct the complement of an embedding. Roughly speaking,
this operation may be done by using homology surgery theory (see [CS] and [V1]) with an obstruction in the relative surgery obstruction group 
$L_*(\Z[\pi_1(\widehat X)]\fl\Lambda)$ where $\Lambda$ is the Cohn localization of the morphism $\Z[\pi_1(\widehat X)]\fl\Z[\pi_1(M)]$.
\vskip 12pt

Let $(X,f)$ be a pointed $B$-space. Suppose $X$ is connected with fundamental group $G$ and $G\fl\pi$ is onto. We have
relative homotopy groups $\pi_*(f)$ appearing in the long exact sequence:
$$\dots\fl \pi_i(X)\fl \pi_i(B)\fl \pi_i(f)\fl \pi_{i-1}(X)\fl\dots$$
The group $\pi_i(f)$ is a $G$-module for $i>2$ and a $G$-group for $i=2$.

Denote by Mod$_G$ the category of right $G$-modules (or equivalently of right $\Z[G]$-modules) and by Gr$_G$ the category of right $G$-groups. An object of Gr$_G$ is a
group equipped with a right $G$-action (compatible with the group structure). In this category we have free objects: a free $G$-group generated by a set $S$ is
the free group generated by $S\times G$.

We have an abelianization functor: Gr$_G\fl$ Mod$_G$ sending finitely generated free $G$-groups to finitely generated free $G$-modules. 

Suppose $\Z[G]\fl R$ is a ring homomorphism. Denote by $W(\Z[G],R)$ the class of morphisms $\varphi:E\fl F$ in Mod$_G$ between two finitely generated free modules sent
to isomorphisms in Mod$_R$ and by $W(G,R)$ the class of morphisms $\varphi:E\fl F$ in Gr$_G$ between two finitely generated free groups sent to isomorphisms in Mod$_R$.
\vskip 12pt
\noi{\bf Theorem 2:} {\sl Let $B$ be a path-connected pointed space with fundamental group $\pi$, $R$ be a ring and $\Z[\pi]\fl R$ be an epimorphism of rings. Let
  $(X,f)$ be a pointed $B$-space and $G$ be the fundamental group of $X$. Suppose $X$ is connected and $G\fl\pi$ is onto. Then $(X,f)$ is $W(B,R)$-local if and only
  if the following holds:

  $\bullet$ $f:X\fl B$ is a Serre fibration

  $\bullet$ $\pi_2(f)$ is a $W(G,R)$-local $G$-group

  $\bullet$ for every integer $i>2$, $\pi_i(f)$ is a $W(\Z[G],R)$-local $G$-module.}
\vskip 12pt
Let $A\fl R$ be a ring homomorphism and $W(A,R)$ be the class of morphisms between finitely generated free $A$-modules sent to isomorphisms in
Mod$_R$. By inverting all morphisms in $W(A,R)$, we get a ring $\Lambda$ called the Cohn localization of $(A\fl R)$ (see [C1] and [V2] section 1). We have a
commutative diagram of rings:
$$\diagram{A&\hfl{}{}&\Lambda\cr\vfl{}{}&&\vfl{}{}\cr A&\hfl{}{}&R\cr}$$
The product map: $\Lambda\build\otimes_A^{}\Lambda\fl\Lambda$ is an isomorphism and we have: Tor$^A_1(M,\Lambda)=0$,
for every right $\Lambda$-module $M$. Moreover, an $A$-module is $W(A,R)$-local if and only if the morphism $M\fl M\build\otimes_A^{}\Lambda$ is an isomorphism and
$(-\build\otimes_A^{}\Lambda,i)$ (where $i$ is induced by to morphism $A\fl\Lambda$) is a localization in Mod$_A$ with respect to $W(A,R)$. This localization will be
denoted by $(L^A_R,i)$.

\vskip 12pt
\noi{\bf Remark:} Let $R'$ be the image of $\Lambda\fl R$. By construction, the three classes $W(A,R)$, $W(A,\Lambda)$ and $W(A,R')$ are the same. If $A$ is the ring
$\Z[G]$, we have also: $W(G,R)=W(G,\Lambda)=W(G,R')$.

Then, up to replacing $R$ by $R'$, we may often suppose that $\Lambda\fl R$ is onto.

\vskip 12pt
\noi{\bf Theorem 3:} {\sl Let $G$ be a group and $\Z[G]\fl R$ be a ring homomorphism.  Then the localization $(L^{\Z[G]}_R,i)$ in Mod$_G$ with respect to $W(\Z[G],R)$ 
can be extended to a localization $(L^G_R,i)$ in Gr$_G$ with respect to $W(G,R)$.}
\vskip 12pt
\noi{\bf Theorem 4:} {\sl Let $B$ be a path-connected pointed space with fundamental group $\pi$, $R$ be a ring and $\Z[\pi]\fl R$ be an epimorphism of rings.

  Let $(X,f)$ be a pointed $B$-space with fundamental group $G$, $(X,f)\fl(Y,g)$ be a $W(B,R)$-localization and $n>1$ be an integer. Suppose $X$ is connected, 
$G\fl\pi$ is onto and $\pi_i(f)$ is $W(G,R)$-local (or $W(\Z[G],R)$-local) for each $i=2,3,\dots,n-1$. Then the morphism $\pi_i(f)\fl\pi_i(g)$ is a bijection for each 
$i<n$ and a $W(G,R)$-localization (or a $W(\Z[G],R)$-localization) for $i=n$.}
\vskip 24pt
\noi{\bf 1 Finite localization of $B$-spaces.}
\vskip 24pt
Throughout this paper $B$ will denote a pointed path-connected space with fundamental group $\pi$ and $\Z[\pi]\fl R$ a ring homomorphism. 

The forgetful functor from the category of $B$-spaces to the category of topological spaces will be denoted by $X\mapsto\underline X$. Then every $B$-space $X$ is a pair 
$X=(\underline X,f)$ where $\underline X$ is a topological space having the homotopy type of a CW complex and $f$ is a continuous map from $\underline X$ to $B$.

The following notations will also be used:

If $X=(\underline X,f)$ is a pointed $B$-space (i.e. $\underline X$ is pointed and $f$ respect basepoints), the fundamental group of $\underline X$ will be denoted by
$\pi_1(X)$ and, for each integer $i>1$, the $i$-th relative homotopy group of $\pi_i(f)$ will be denoted by $\Gamma_i(X)$. Then $\Gamma_i(X)$ is a $\pi_1(X)$-group (i.e. 
a group with a right $G$-action) for $i=2$ and a $\pi_1(X)$-module for $i>2$. We have also an exact sequence of $\pi_1(X)$-groups:
$$\pi_2(B)\fl\Gamma_2(X)\fl\pi_1(X)\fl \pi$$

We say that a $B$-space $X=(\underline X,f)$ is a $B$-complex or a finite $B$-complex if $\underline X$ is a CW complex or a finite CW complex. We say also that
a map $X\fl Y$ is an inclusion of $B$-complexes if $X$ and $Y$ are $B$-complexes and the induced map $\underline X\fl \underline Y$ is an inclusion of CW complexes.

The class of inclusions $X\fl Y$ of finite $B$-complexes such that $\H_*(\underline Y,\underline X;R)=0$ will be denoted by $W(B,R)$.

If $W$ is a class of cofibrations in Sp$_B$, we say that a $B$-space is weakly $W$-local if, for every morphism $Y\fl Z$ in $W$, every morphism $Y\fl X$ has an
extension to $Z$.

\vskip 12pt
\noi{\bf 1.1 Lemma:} {\sl Let $X$ be a weakly $W(B,R)$-local $B$-space. Then the map $f:\underline X\fl B$ is a Serre fibration and $X$ is $W(B,R)$-local.}
\vskip 12pt
\noi{\bf Proof:} Let $n\geq0$ be an integer. Consider a commutative diagram:
$$\diagram{B^n&\hfl{i_0}{}&B^n\times[0,1]\cr \vfl{g}{}&&\vfl{h}{}\cr \underline X&\hfl{f}{}&B\cr}$$
where $B^n$ is the $n$-ball and $i_0$ is the map $u\mapsto(u,0)$. This map is a map from $(B^n,fg)$ to $(B^n\times[0,1],h)$ and belongs clearly to $W(B,R)$. Since $X$ is 
weakly $W(B,R)$-local, there is a map $\varphi:B^n\times[0,1]\fl \underline X$ such that:
$$\varphi i_0=g\hskip 24pt\hbox{and}\hskip 24pt f\varphi=h$$
and $f$ is a Serre fibration.

Let $\alpha:K\fl L$ be a morphism in $W(B,R)$. This map induces a map $\alpha^*:\Map(L,X)\fl\Map(K,X)$. 
Since $\alpha$ is a cofibration and $\underline X\fl B$ a Serre fibration, $\alpha^*$ is also a Serre fibration. Then,
in order to prove that $X$ is $W(B,R)$-local, the only thing to do is to prove that each fiber of $\alpha^*$ is weakly contractible.

Let $g:K\fl X$ be a map and $F_g$ be the corresponding fiber of $\alpha^*$. Since $X$ is weakly $W(B,R)$-local, the map $\alpha^*$ is epic and $F_g$ is nonempty.

Let $n\geq0$ be an integer and $u:S^n\fl F_g$ be a continuous map. So we have a commutative diagram:
$$\diagram{S^n&\hfl{u}{}&\Map(L,X)\cr\vfl{}{}&&\vfl{\alpha^*}{}\cr pt&\hfl{g}{}&\Map(K,X)\cr}$$
which is equivalent to the following commutative diagram of $B$-spaces:
$$\diagram{S^n\times K&\hfl{}{}&B^{n+1}\times K\cr\vfl{Id\times\alpha}{}&&\vfl{g\circ pr_2}{}\cr S^n\times L&\hfl{v}{}&X\cr}$$
Denote by $E$ the subspace $B^{n+1}\times K\cup S^n\times L$ of $B^{n+1}\times L$. This space is a $B$-space defined by the cocartesian diagram:
$$\diagram{S^n\times K&\hfl{}{}&B^{n+1}\times K\cr\vfl{}{}&&\vfl{i}{}\cr S^n\times L&\hfl{j}{}&E\cr}$$
and there is a unique morphism $w:E\fl X$ such that: $wi=g\circ pr_2$ and $wj=v$.

On the other hand, the inclusion $E\subset B^{n+1}\times L$ belongs to the class $W(B,R)$ and $w$ has an extension $w':B^{n+1}\times L\fl X$.
This map corresponds to a continuous map $w'':B^{n+1}\fl\Map(L,X)$. The two conditions $wi=g\circ pr_2$ and $wj=v$ imply that $w''$ is a map $B^{n+1}\fl F_g$ which is
an extension of $u:S^n\fl F_g$.

Hence every map $S^n\fl F_g$ is homotopic to a point and $F_g$ is weakly contractible. Since this property holds for each fiber of $\alpha^*$,
the $B$-space $X$ is $W(B,R)$-local.\cqfd
\vskip 12pt
Let $\Omega$ be a set of finite CW-complexes such that every finite CW-complex is isomorphic to a complex in $\Omega$.

If $X$ is a $B$-space, denote by $V(X)$ the set of diagrams in Sp$_B$:
$$\diagram{K&\hfl{\alpha}{}&L\cr\vfl{\beta}{}&&\cr X&&\cr}$$
where $\underline K$ and $\underline L$ are in $\Omega$ and $\alpha$ in $W(B,R)$.

By taking the coproduct of all these morphisms $K\fl L$, we get a cofibration $\coprod K\fl\coprod L$ and a cocartesian diagram:
$$\diagram{\coprod K&\hfl{}{}&\coprod L\cr\vfl{}{}&&\vfl{}{}\cr X&\hfl{j}{}&E(X)\cr}$$

Since $X\mapsto V(X)$ is clearly a functor from Sp$_B$ to the category of sets, this construction is functorial. In particular $E$ is a functor from Sp$_B$ to
itself and $j:$ Id$\fl E$ is a morphism of functors. Moreover, for every $B$-space $X$, $j:X\fl E(X)$ is a cofibration.

Using this functor, we get a sequence of $B$-spaces:
$$X=X_0\build\fl_{}^j X_1\build\fl_{}^j X_2\build\fl_{}^j X_3\build\fl_{}^j \dots$$
where, for every $n\geq0$, $X_{n+1}$ is the $B$-space $E(X_n)$. Moreover all these maps $j$ are cofibrations and this sequence has a colimit denoted
by $L^B_R(X)$. So we get a functor $L^B_R$ and a morphism $i:$ Id$\fl L^B_R$ such that $i:X\fl L^B_R(X)$ is a cofibration for every $B$-space $X$. 
\vskip 12pt
\noi{\bf 1.2 Lemma:} {\sl For every $B$-space $X$, $L^B_R(X)$ is $W(B,R)$-local.}
\vskip 12pt
\noi{\bf Proof:} Consider a diagram in Sp$_B$:
$$\diagram{K&\hfl{\alpha}{}&L\cr\vfl{\beta}{}&&\cr L^B_R(X)&&\cr}$$
where $\alpha$ is in $W(B,R)$. Since $\underline K$ is a finite CW complex, there is an integer $n$ such that $\beta(K)\subset X_n$ and we have a diagram:
$$\diagram{K&\hfl{\alpha}{}&L\cr\vfl{\beta}{}&&\cr X_n&&\cr}$$
This diagram is isomorphic to some element in $V(X_n)$. Therefore the map $K\fl X_n\subset X_{n+1}$ factors through $L$ and the map $K\build\fl_{}^\beta L^B_R(X)$
has an extension to $L$.

Thus $L^B_R(X)$ is weakly $W(B,R)$-local and, because of lemma 1.1, $L^B_R(X)$ is $W(B,R)$-local.\cqfd
\vskip 12pt
\noi{\bf 1.3 Lemma:} {\sl For every $B$-space $X$ the map $i:X\fl L^B_R(X)$ is a $W(B,R)$-equivalence.}
\vskip 12pt
\noi{\bf Proof:} Let $Y$ be a $W(B,R)$-local $B$-space. Then we have a sequence of maps:
$$\Map(X_0,Y)\longleftarrow\Map(X_1,Y)\longleftarrow\Map(X_2,Y)\longleftarrow\Map(X_3,Y)\longleftarrow\dots$$
and, because each morphism $X_n\fl X_{n+1}$ is a cofibration, each restriction map $\Map(X_j,Y)\fl\Map(X_i,Y)$ is a Serre fibration.

Let $n\geq0$ be an integer and $g$ be a map from $X_n$ to $Y$. Denote by $F(g)$ the fiber of $g$ in the restriction map: $\Map(X_{n+1},Y)\fl\Map(X_n,Y)$.

By construction, there is a family of morphisms $\alpha_i:K_i\fl L_i$ in $W(B,R)$ and morphisms $\beta_i:K_i\fl X_n$ such that $X_{n+1}$ is obtained by the cocartesian
diagram:
$$\diagram{K&\hfl{\alpha}{}&L\cr\vfl{\beta}{}&&\vfl{}{}\cr X_n&\hfl{}{}&X_{n+1}\cr}$$
where $\alpha:K\fl L$ is the coproduct of the $\alpha_i$'s and $\beta:K\fl X_n$ is induced by the $\beta_i$'s. Therefore the fiber $F(g)$ is the product of the spaces
$F(\beta_i\circ g)$, where $F(\beta_i\circ g)$ is the fiber of $\beta_i\circ g$ in the restriction map $\alpha_i^*:\Map(L_i,Y)\fl\Map(K_i,Y)$. Since $\alpha_i$ belongs
to the class $W(B,R)$ each fiber $F(\beta_i\circ g)$ is weakly contractible and so is $F(g)$.

Hence each fiber of $\Map(X_{n+1},Y)\fl\Map(X_n,Y)$ is weakly contractible and so is each fiber of $\Map(X_n,Y)\fl\Map(X,Y)$. By passing to the limit we see that
$\Map(L^B_R(X),Y)\fl\Map(X,Y)$ is a Serre fibration with weakly contractible fibers. The result follows.\cqfd
\vskip 12pt
These two lemmas imply that $(L^B_R,i)$ is a localization in Sp$_B$ with respect to the class $W(B,R)$. Therefore a $B$-space $X$ is
$W(B,R)$-local if and only if $i:X\fl L^B_R(X)$ is a strong deformation retract.
\vskip 12pt
\noi{\bf 1.4 Lemma:} {\sl Let $g$ be a map from a $B$-complex $X$ to a $B$-space $Y$. Then there exists a commutative diagram:
  $$\diagram{X&\hfl{}{}&Z\cr\vfl{=}{}&&\vfl{}{}\cr X&\hfl{g}{}&Y\cr}$$
  where $Z$ is a $B$-complex, $X\fl Z$ an inclusion and $Z\fl Y$ a homotopy equivalence.}
\vskip 12pt
\noi{\bf Proof:} Since $Y$ is a $B$-space, there is a CW complex $H$ and a homotopy equivalence $h:H\fl \underline Y$. Since $\underline X$ is a
CW complex, there is a cellular map $\alpha:\underline X\fl H$ such that $h\alpha$ is homotopic to $g:\underline X\fl\underline Y$. Denote by $C$ the
cylinder of the map $\alpha$. The homotopy between $h\alpha$ and $g$ induces a map $C\fl\underline Y$ extending $h$ and we have a commutative diagram:
$$\diagram{\underline X&\hfl{}{}&C\cr\vfl{=}{}&&\vfl{}{}\cr\underline X&\hfl{}{}&\underline Y\cr}$$
where $C$ is a CW complex, $\underline X\fl C$ is an inclusion of CW complexes and $C\fl\underline Y$ is a homotopy equivalence. Therefore the diagram above can be
lifted in the desired diagram of $B$-spaces where $Z$ is a $B$-complex with underline space $C$.\cqfd
\vskip 12pt
\noi{\bf 1.5 Lemma:} {\sl For every $B$-space $X$, the morphism $j:X\fl E(X)$ belongs to the class $\overline W(B,R)$.}
\vskip 12pt
\noi{\bf Proof:} Let $K\subset L$ be an inclusion of finite $B$-complexes and $\varphi:(K\subset L)\fl (X\build\fl_{}^j E(X))$ be a morphism in Mor$($Sp$_B)$.

By construction, $E(X)$ is defined by a cocartesian diagram:
$$\diagram{U&\hfl{\alpha}{}&V\cr\vfl{\beta}{}&&\vfl{}{}\cr X&\hfl{j}{}&E(X)\cr}$$
where $\alpha:U\fl V$ is a disjoint union of morphisms $\alpha_i:U_i\fl V_i$ in $W(B,R)$.

By applying lemma 1.4 to the map $K\coprod U\fl X$ induced by $\varphi$ and $\beta$, we get a $B$-complex $Y$, two inclusions $\varphi_0:K\fl Y$ and $\gamma:U\fl Y$
and a homotopy equivalence $g:Y\fl X$ such that: $g\gamma=\beta$ and $g\varphi_0$ is the map $\varphi:K\fl X$.

Let $Z$ be the push-out of $\alpha:U\fl V$ and $\gamma:U\fl Y$. We have two cocartesian diagrams:
$$\diagram{U&\hfl{}{}&V\cr\vfl{}{}&&\vfl{}{}\cr Y&\hfl{}{}&Z\cr}\hskip 48pt\diagram{Y&\hfl{}{}&Z\cr\vfl{g}{}&&\vfl{h}{}\cr X&\hfl{}{}&E(X)\cr}$$
Since $g:Y\fl X$ is a homotopy equivalence, the induced map $\psi:(Y\fl Z)\fl (X\fl E(X))$ induced by the second diagram is also a homotopy equivalence.

Since $Z\fl E(X)$ is a homotopy equivalence, there is a cellular map $\delta:\underline L\fl\underline Z$ such that $h\delta:\underline L\fl\underline{E(X)}$ is
homotopic to $\varphi:\underline L\fl\underline{E(X)}$. So we get a commutative diagram:
$$\diagram{\underline L&\hfl{\lambda_0}{}&\underline L\times[0,1]&\mfl{\lambda_1}{}&\underline L\cr
  \vfl{\varphi}{}&&\vfl{\psi}{}&&\vfl{\delta}{}\cr \underline{E(X)}&\hfl{=}{}&\underline{E(X)}&\mfl{h}{}&\underline Z\cr}$$
where $\lambda_t$ is the map $x\mapsto(x,t)$. This diagram induces a commutative diagram of $B$-spaces:
$$\diagram{ L&\hfl{\lambda_0}{}&L_1&\mfl{\lambda_1}{}&L\cr
  \vfl{\varphi}{}&&\vfl{\psi}{}&&\vfl{\delta}{}\cr E(X)&\hfl{=}{}&E(X)&\mfl{h}{}&Z\cr}$$
where $K_1$ is a $B$-complex such that: $\underline L_1=L\times[0,1]$.

Let $J$ be a finite subset of the set $I$ of indices $i$ and $F$ be a finite $B$-subcomplex of $Y$ containing $\beta(U_j)$, for all $j\in J$. We
have a cocartesian diagram:
$$\diagram{U(J)&\hfl{}{}&V(J)\cr\vfl{}{}&&\vfl{}{}\cr F&\hfl{}{}&F'\cr}$$
where $U(J)$ and $V(J)$ are the unions of the $U_j$ and $V_j$, for all $j\in J$. Moreover, $F'$ is a finite $B$-subcomplex of $Z$ and the inclusion $F\subset F'$
belongs to the class $W(B,R)$. By taking $J$ and $F$ big enough, we have: $\varphi_0(K)\subset F$ and: $\delta(L)\subset F'$. Let $F'_1$ be the push-out of
$\lambda_1:L\fl L_1$ and $\delta:L\fl F'$.

We have a commutative diagram:
$$\diagram{K&&L&\hfl{\lambda_1}{}&L_1\cr \vfl{\varphi_0}{}&&\vfl{\delta}{}&&\vfl{}{}\cr F&\hfl{}{}&F'&\hfl{}{}&F'_1\cr\vfl{g}{}&&\vfl{h}{}&&\vfl{}{}\cr
  X&\hfl{}{}&E(X)&\hfl{=}{}&E(X)\cr}$$

Since $\lambda_1$ is a homotopy equivalence, the map $F'\fl F'_1$ is also a homotopy equivalence and the morphism $F\fl F'_1$ belongs to the class $W(B,R)$.

The morphism $\varphi:(K\subset L)\fl(X\fl E(X))$ factors through $(K\subset L_1)$ (via the map $\lambda_0:L\fl L_1$) and then through $(F\fl F'_1)$ which is in
$W(B,R)$. Therefore $j:X\fl E(X)$ belongs to the class $\overline W(B,R)$.\cqfd
\vskip 12pt
\noi{\bf 1.6 Lemma:} {\sl The class $\overline W(B,R)$ is stable under composition.}
\vskip 12pt
\noi{\bf Proof:} Let $X\fl Y$ and $Y\fl Z$ be two morphisms in $\overline W(B,R)$. Let $(K\subset Y)$ be an inclusion of finite $B$-complexes and
$\varphi:(K\subset L)\fl(X\fl Z)$ be a morphism in Mor$($Sp$_B)$.

This morphism $\varphi$ induces a morphism $(K\subset L)\fl(Y\fl Z)$ which factors through a morphism $(H_1\subset L_1)$ in $W(B,R)$. So we get a morphism
$(K,H_1)\fl (X\subset Y)$ which factors through a morphism $(K_2\subset H_2)$ in $W(B,R)$.

By applying lemma 1.4 to the map $K_2\coprod H_1\fl H_2$, we get a $B$-complex $H'_2$, a homotopy equivalence $H'_2\fl H_2$ and inclusions $K_2\subset H'_2$ and
$H_1\subset H'_2$. Let $L_2$ be the push-out of $H_1\fl H'_2$ and $H_1\fl L_1$. Then we have a commutative diagram:
$$\diagram{K&\hfl{}{}&H_1&\hfl{}{}&L_1\cr\vfl{}{}&&\vfl{}{}&&\vfl{}{}\cr K_2&\hfl{}{}&H'_2&\hfl{}{}&L_2\cr\vfl{}{}&&\vfl{}{}&&\vfl{}{}\cr X&\hfl{}{}&Y&\hfl{}{}&Z\cr}$$

Since $K_2\fl H'_2$ and $H'_2\fl L_2$ are in $W(B,R)$, so is the same for $K_2\fl L_2$. Hence the map $\varphi:(K\subset L)\fl(X\fl Z)$ factors through $(K\subset L_1)$
and then through $(K_2\fl L_2)$ which is $W(B,R)$. Therefore $X\fl Z$ belongs to the class $\overline W(B,R)$ and the lemma is proven.\cqfd
\vskip 12pt
\noi{\bf 1.7 Lemma:} {\sl For every $B$-space $X$, the inclusion $i:X\fl L^B_R(X)$ belongs to the class $\overline W(B,R)$.}
\vskip 12pt
\noi{\bf Proof:} $L^B_R(X)$ is the colimit of the sequence: $X=X_0\subset X_1\subset X_2\subset X_3\subset\dots$, with $X_{n+1}=E(X_n)$. Because of lemmas 1.5 and
1.6, the inclusion $X\subset X_n$ belongs to $\overline W(B,R)$ for every $n$.

Let $(K\subset L)$ be an inclusion of finite $B$-complexes and $\varphi:(K\subset L)\fl(X\fl L^B_R(X))$ be a map. Since $L$ is finite $\varphi$ factors through some
$(X\subset X_n)$ and then through a morphism in $W(B,R)$. The result follows and theorem 1 is proven.\cqfd
\vskip 12pt
\noi{\bf 1.8 Remark:} Let $A\fl R$ be a homomorphism of rings. Denote by $\C_*(A)$ the category of $A$-complexes (i.e. the category of graded differential projective
right $A$-modules). An $A$-complex $C$ is connective if $C_i=0$ for all $i<0$. The category $\C'_*(A)$ of connective $A$-complexes is a subcategory of $\C_*(A)$.

In these categories we have notions of cofibration, homotopy and finite objects: a $A$-complex $C$ is said to be finite if each $C_i$ is free and $\oplus_i C_i$ is
finitely generated.

Using the morphism $A\fl R$, we have also two classes of morphisms $W_*(A,R)\subset \C_*(A)$ and $W'_*(A,R)\subset \C'_*(A)$: a morphism $\varphi:C\fl C'$ in $\C(A)$
(or in $\C'(A)$) belongs to $W_*(A,R)$ (or in $W'_*(A,R)$) if $C$ and $C'$ are finite and $\varphi$ induces an isomorphism $\H_*(C,R)\build\fl_{}^\sim\H_*(C',R)$.

Therefore, by using exactly the same method as above, it is possible to construct two localization functors $(L^A_{*R},i)$ and $(L'^A_{*R},i)$ in $\C_*(A)$ and
$\C'_*(A)$ with respect to the classes $W_*(A,R)$ and $W'_*(A,R)$.
\vskip 24pt
\noi{\bf 2 Local $B$-spaces and homotopy groups.}
\vskip 12pt
This section is devoted to the proof of theorem 2.

Let $X=(\underline X,f)$ be a pointed $B$-space. Denote by $G=\pi_1(X)$ the fundamental group of $\underline X$ and suppose $\underline X$ is connected and $G\fl\pi$ is 
onto. We have to prove that $X$ is $W(B,R)$-local if and only if $f$ is a Serre fibration, $\Gamma_2(X)=\pi_2(f)$ is $W(G,R)$-local and $\Gamma_i(X)=\pi_i(f)$ is 
$W(\Z[G],R)$-local for every $i>2$.
\vskip 12pt
\noi{\bf 2.1 Lemma:} {\sl Suppose $X$ is $W(B,R)$-local. Then $f:\underline X\fl B$ is a Serre fibration, $\Gamma_2(X)$ is $W(G,R)$-local and $\Gamma_i(X)$ is 
$W(\Z[G],R)$-local for every integer $i>2$.}
\vskip 12pt
\noi{\bf Proof:} Because of lemma 1.1, $X$ is weakly $W(B,R)$-local and $f$ is a Serre fibration.

Let $k>1$ be an integer. Denote by $W$ the class $W(G,R)$ is $k=2$ and the class $W(\Z[G],R)$ if $k>2$.

If $G'\fl G$ is a group homomorphism and $S$ is a set, we let $F_{G'}(S)$ denote the $G'$-group freely generated by $S$ if $k=2$ and the $Z[G']$-module freely generated by
$S$ if $k>2$.

Consider two finite sets $U$ and $V$ and a diagram of $G$-groups or $\Z[G]$-modules:
$$\diagram{F_G(U)&\hfl{\alpha}{}&F_G(V)\cr\vfl{\lambda}{}&&\cr\Gamma_k(X)&&\cr}$$
with $\alpha$ in $W$. We have to prove that $\lambda:F_G(U)\fl\Gamma_k(X)$ has a unique extension $\mu:F_G(V)\fl\Gamma_k(X)$.

Consider a pointed $B$-complex $K$ and a pointed map $K\fl X$. Suppose $K$ is finite and connected. So we have two morphisms $\Gamma_k(K)\fl\Gamma_k(X)$ 
and $F_{G'}(U)\fl F_G(V)$ with $G'=\pi_1(K)=\pi_1(\underline K)$. Up to taking $K$ big enough, we may as well suppose that $\lambda(U)$ is contained in the image of  
$\Gamma_k(K)\fl\Gamma_k(X)$ and that $\alpha(U)$ is contained in the image of $F_{G'}(U)\fl F_G(V)$. Hence we may lift $\lambda:U\fl\Gamma_k(X)$ to a map 
$\lambda':U\fl \Gamma_k(K)$ and $\alpha:U\fl F_G(V)$ to a map $\alpha':U\fl F_{G'}(V)$. So we get two morphisms $\lambda':F_{G'}(U)\fl \Gamma_k(K)$ and 
$\alpha':F_{G'}(U)\fl F_{G'}(V)$ and two commutative diagrams:
$$\diagram{F_{G'}(U)&\hfl{}{}&F_G(U)\cr\vfl{\lambda'}{}&&\vfl{\lambda}{}\cr \Gamma_k(K)&\hfl{}{}&\Gamma_k(X)\cr}\hskip 48pt
\diagram{F_{G'}(U)&\hfl{\alpha'}{}&F_{G'}(V)\cr\vfl{}{}&&\vfl{}{}\cr F_G(U)&\hfl{\alpha}{}&F_G(V)\cr}$$

Denote by $S(V)$ a bouquet of pointed spheres $S^{k-1}$ indexed by $V$. With the trivial map $S(V)\fl B$, $S(V)$ becomes a $B$-complex and we get a $B$-complex
$K_1=K\vee S(V)$. 

For each $v\in V$, the corresponding sphere will be denoted by $S_v$ and the map $S_v\subset S(V)\subset K_1$ defines an element in the $k-1$ homotopy group of the 
homotopy fiber of $\underline K_1\fl B$ and then an element $\mu(v)\in\Gamma_k(K_1)$. This correspondence extends to a morphism of $G'$-groups (or $\Z[G']$-modules) 
$\mu:F_{G'}(V)\fl\Gamma_k(K_1)$.

For each $u\in U$, we have an element $\beta(u)\in\pi_k(g_1)$ defined by: $\beta(u)=\mu(\alpha'(u))\lambda'(u)^{-1}$ if $k=2$ and by: 
$\beta(u)=\mu(\alpha'(u))-\lambda'(u)$ if $k>2$. For each $u\in U$, we may add a $k$-cell $B_u$ to $K_1$ with the attaching map $\beta(u)$. So we get a $B$-complex $L$ and
a diagram of $B$-spaces:
$$\diagram{K&\hfl{}{}&L\cr\vfl{}{}&&\cr X&&\cr}$$

Moreover the cellular complex of $(K\subset L)$ reduces to the morphism $F_{G'}(U)\build\fl_{}^{\alpha'}F_{G'}(V)$ in degree $k-1$ and $k$. But this morphism induces
the morphism $F_G(U)\build\fl_{}^{\alpha}F_G(V)$ which is in the class $W$. Therefore the map $K\fl L$ belongs to the class $W(B,R)$ and, because $X$ is $W(B,R)$-local,
there exists a map $h:L\fl X$, unique up to homotopy, extending the map $K\fl X$.

Such an extension is defined by an extension $h_*:\Gamma_k(K_1)\fl\Gamma_k(X)$ of $\Gamma_k(K)\fl\Gamma_k(X)$ sending $\beta(u)$ to $1$ or $0$, for each $u\in U$. This 
extension is determined by a map $\nu:U\fl\Gamma_k(X)$ sending each $u\in U$ to $\nu(u)=h_*(\mu(u))$. Therefore there is a unique morphism $\nu':F_{G'}(U)\fl\Gamma_k(X)$ 
such that:
$$\nu'(\alpha'(u))\lambda'(u)^{-1}=1\hskip 24pt\hbox{or}\hskip 24pt\nu'(\alpha'(u))-\lambda'(u)=0$$
for each $u\in U$. Hence there is a unique morphism $\nu':F_{G'}(U)\fl\Gamma_k(X)$ such that: $\nu'\alpha'=\lambda'$ and then a unique morphism $\nu:F_G(U)\fl\Gamma_k(X)$ 
such that $\nu\alpha=\lambda$. Thus each $\Gamma_k(X)$ is $W$-local and the lemma is proven.\cqfd
\vskip 12pt
Suppose now $f$ is a Serre fibration, $\Gamma_2(X)$ is $W(G,R)$-local and $\Gamma_i(X)$ is $W(\Z[G],R)$-local for every $i>2$. We have to prove that $X$ is $W(B,R)$-local.
In order to do that, we will need other classes of morphisms in Sp$_B$:

Let $\widehat W$ be the class of morphisms $K\fl L$ in Sp$_B$ with the following properties:

$\bullet$ $L$ is a $B$-complex, $K$ is a subcomplex of $L$ and $\underline L/\underline K$ is finite

$\bullet$ the inclusion $K\fl L$ induces a bijection $\pi_0(\underline K)\simeq\pi_0(\underline L)$ and an isomorphism $\H_*(K,R)\build\fl_{}^\sim\H_*(L,R)$.

If $p<q$ are two integers, we'll denote by $W(p,q)$ the class of morphisms $K\fl L$ in $\widehat W$ such that the dimension of each cell in $L\setminus K$ is in 
$\{p,p+1,\dots,q\}$. For simplicity, the class $W(n-1,n)$ will also be denoted by $W(n)$.

Let $K\fl L$ be a morphism in $W(B,R)$ and $g:K\fl X$ be a morphism in Sp$_B$. We have to prove that $g:K\fl X$ has an extension $L\fl X$.

The morphisms $K\fl L$ and $K\fl X$ induce a morphism $\varphi$ in Mor$($Sp$_B)$ from $K\fl L$ to $X\fl B$, where $B$ is considered as the terminal object in Sp$_B$.

We will construct $B$-complexes $K_1$, $K_2$, \dots, $K_n$, $L'$, morphisms:
$$K_1\fl K_2\fl \dots \fl K_n\fl L'$$
and a factorization of $\varphi$ through $K_1\fl L'$ with the following properties:

$\bullet$ $K_1$, $K_2$, \dots, $K_n$ and $L'$ are connected pointed $B$-complexes and the morphisms $K_i\fl K_{i+1}$, $K_n\fl L'$ and $K_1\fl X$ are pointed maps

$\bullet$ the morphism $K_1\fl X$ induces bijections $\pi_1(K_1)\build\fl_{}^\sim G$ and $\Gamma_2(K_1)\build\fl_{}^\sim \Gamma_2(X)$

$\bullet$ the morphism $K_1\fl L'$ belongs to $W(1,n)$, the morphism $K_i\fl K_{i+1}$ belongs, for each $i<n$, to the class $W(i+1)$ and $K_n\fl L'$ is a strong homotopy 
equivalence.

After that we'll prove that the map $K_1\fl X$ has an extension to $K_2$ and then to $K_n$ and to $L'$. This last extension induces the desired extension of $g:K\fl X$ to
the $B$-complex $L$. 
\vskip 12pt
\noi{\bf 2.2 Construction of $K_1$ and $L'$.}
\vskip 12pt
Since $K\fl L$ induces a bijection $\pi_0(\underline K)\build\fl_{}^\sim\pi_0(\underline L)$, we may add to $\underline K$ $0$ and $1$-cells in $\underline L$ in order to
get a CW-complex $\underline K'$ contained in $\underline L$ with the following properties:

$\bullet$ $\underline K\fl\underline K'$ is a homotopy equivalence

$\bullet$ every vertex in $\underline L$ is in $\underline K'$.

Since $f$ is a Serre fibration, the map $\underline K\fl\underline X$ has an extension to $\underline K'$ and we get a $B$-complex $K'$ such that: $K\subset K'\subset L$.
Moreover the morphism $\varphi$ from $K\fl L$ to $X\fl B$ factors through $K'\fl L$ and there is an integer $n$ such that $K'\fl L$ belongs to the class $W(1,n)$.

By adding cells to $\underline K'$ we may get a pointed CW complex $\underline K_1$ and an extension $\underline K_1\fl \underline X$ of $\underline K'\fl\underline X$ 
such that the map $\underline K_1\fl \underline X$ is pointed and $2$-connected. So we get a pointed $B$-complex $K_1$ containing $K'$ and a factorization of $K'\fl X$ 
through $K_1$. Moreover $K_1$ is connected and the map $K_1\fl X$ induces bijections $\pi_1(K_1)\build\fl_{}^\sim G$ and $\Gamma_2(K_1)\build\fl_{}^\sim \Gamma_2(X)$.

The $B$-complex $L'$ is defined by the cocartesian diagram:
$$\diagram{K'&\hfl{}{}&L\cr\vfl{}{}&&\vfl{}{}\cr K_1&\hfl{\alpha_1}{}&L'\cr}$$
So the map $\varphi$ from $K\fl L$ to $X\fl B$ factors through $K'\fl L$ and then through $K_1\fl L'$. Moreover the inclusion $\alpha_1:K_1\fl L'$ belongs to $W(1,n)$.
\vskip 12pt
\noi{\bf 2.3 Construction of the $K_i$'s.}
\vskip 12pt
The $K_i$'s will be constructed by induction. Let $p$ be an integer with $0<p<n$. Suppose constructed $B$-complexes $K_1$, $K_2$,\dots, $K_p$ and $L_p$, inclusions
$K_1\subset K_2\subset\dots\subset K_p\subset L_p$ and a strong homotopy equivalence $L_p\fl L'$ with the following properties:

$\bullet$ for each $i\in\{2,3,\dots,p\}$, the inclusion $K_{i-1}\subset K_i$ belongs to $W(i)$, 

$\bullet$ the inclusion $K_p\subset L_p$ belongs to $W(p,n)$

$\bullet$ $\alpha_1:K_1\fl L'$ is the composite map $K_1\subset K_2\subset\dots\subset K_p\subset L_p\fl L'$.

These properties are satisfied for $p=1$ with $L_1=L'$.

Suppose $p<n-1$. We have to define $K_{p+1}$ and $L_{p+1}$. We'll denote by $\alpha_p$ the inclusion $K_p\subset L_p$, by $\Sigma$ the set of $p$-cells in 
$L_p\setminus K_p$ and by $\Sigma'$ the set of $p+1$-cells in $L_p\setminus K_p$. Since $K_p\subset L_p$ belongs to $W(p,n)$, $\Sigma$ and $\Sigma'$ are finite sets.

The union of $K_p$ and the cells in $\Sigma$ will be denoted by $K'_p$ and the union of $K'_p$ and the cells in $\Sigma'$ will be denoted by $K''_p$.

Let's take an orientation for each cell in $\Sigma$ and $\Sigma'$. Then each cell $\sigma\in\Sigma'$ induces an element $\lambda(\sigma)$ in 
$\pi_{p+1}(\underline K''_p,\underline K'_p)$.

For every set $S$, $F_G[S]$ will denote the $G$-group freely generated by $S$ and, if $A$ is a ring, the right $A$-module freely generated by $S$ will be denoted by 
$A[S]$.

If $p=1$, $\lambda$ induces a $G$-group homomorphism $\lambda: F_G[\Sigma']\fl\pi_2(\underline K''_1,\underline K'_1)$. If $p>1$, we get a morphism of $\Z[G]$-modules
$\lambda:\Z[G][\Sigma']\fl\pi_{p+1}(\underline K''_p,\underline K'_p)$.

On the other hand, we have a boundary homomorphism:
$$\pi_{p+1}(\underline K''_p,\underline K'_p)\fl\pi_p(\underline K'_p,\underline K_p)\fl\H_p(\underline K'_p,\underline K_p;\Z[\pi])\fl\H_p(\underline K'_p,\underline 
K_p;R)\build\fl_{}^\sim R[\Sigma]$$

According to the value of $p$, we have one of the two commutative diagrams:
$$\diagram{F_G[\Sigma']&\hfl{\lambda}{}&\pi_2(\underline K'_1,\underline K_1)\cr\vfl{}{}&&\vfl{}{}\cr R[\Sigma']&\hfl{d}{}&R[\Sigma]\cr}\hskip 48pt
\diagram{\Z[G][\Sigma']&\hfl{\lambda}{}&\pi_{p+1}(\underline K'_p,\underline K_p)\cr\vfl{}{}&&\vfl{}{}\cr R[\Sigma']&\hfl{d}{}&R[\Sigma]\cr}$$

But $d:R[\Sigma']\fl R[\Sigma]$ is the boundary homomorphism of the cellular complex $C_*(\underline L_p,\underline K_p;R)$ which is acyclic because $K_p\subset L_p$ 
belongs to the class $W(p,n)$. Therefore the morphisms $F_G[\Sigma']\fl\Z[G][\Sigma']\fl R[\Sigma']\fl R[\Sigma]$ are all surjective and every $\sigma\in\Sigma$ may be 
lifted to an element $s(\sigma)$ in $F_G[\Sigma']$ if $p=1$ or in $\Z[G][\Sigma']$ is $p>2$.

Hence we may attached for each $\sigma\in\Sigma$ a $p+1$-cell to $\underline K'_p$ with the map $\lambda s(\sigma)$. So we get a $B$-complex $K_{p+1}$ containing $K_p$ as
a subcomplex. Moreover the elements $s(\sigma)$ induce a well defined map $\mu: K_{p+1}\fl K''_p$. By construction, the cellular complex of 
$(\underline K_{p+1},\underline K_p)$, with coefficients in $R$, is isomorphic to $\dots\fl0\fl R[\Sigma]\build\fl_{}^=R[\Sigma]\fl0\fl\dots$ and the inclusion
$K_p\subset K_{p+1}$ belongs to $W(p+1)$.

Let $L'_{p+1}$ denote the cylinder of the map $\mu:K_{p+1}\fl L_p$. We have a cocartesian square of $B$-spaces:
$$\diagram{K_p\times[0,1]&\hfl{}{}&L'_{p+1}\cr\vfl{pr_1}{}&&\vfl{}{}\cr K_p&\hfl{}{}&L_{p+1}\cr}$$
inducing an inclusion $K_{p+1}\subset L_{p+1}$ of $B$-complexes and a strong homotopy equivalence $L_{p+1}\fl L_p$. Moreover, the transformation $L_p\mapsto L_{p+1}$ is
obtained by deleting the $p$-cells in $L_p\setminus K_p$ and adding $p+1$ and $p+2$-cells. Therefore the inclusion $K_{p+1}\subset L_{p+1}$ belongs to $W(p+1,n)$.

Thus this construction may be iterated and we get $B$-complexes $K_1\subset K_2\subset\dots\subset K_{n-1}\subset L_{n-1}$ and a strong homotopy equivalence 
$L_{n-1}\fl L'$. But the inclusion $K_{n-1}\subset L_{n-1}$ belongs to $W(n-1,n)=W(n)$ and we may set: $K_n=L_{n-1}$.
\vskip 12pt
\noi{\bf 2.4 Construction of the map $g_2:K_2\fl X$.}
\vskip 12pt
Let $\Sigma$ be the set of $1$-cells in $K_2\setminus K_1$. As before, we will take an orientation for each cell in $\Sigma$.

Consider a cell $\sigma\in\Sigma$. This cell is an oriented edge in $K_2$. Since $\pi_1(K_1)\fl\pi$ is onto, there are two paths in $\underline K_1$ joining
the base point in $\underline K_1$ to the end points of $\sigma$ such that the union $\widehat\sigma$ of $\sigma$ and these two paths is an oriented loop in 
$\underline K_2$ which is homotopic to the trivial loop in $B$.

Let $S$ be a bouquet of circles indexed by $\Sigma$. For each $\sigma\in\Sigma$, the corresponding circle in $S$ will be denoted by $S(\sigma)$. Denote by 
$\underline H'$ the space $\underline K_1\vee S$. The identity of $\underline K_1$ and the maps $\widehat\sigma:S(\sigma)\fl K_2$ induce a map $\alpha:\underline H'\fl
\underline K_2$. Denote by $K'$ the union of $K_1$ and the $1$-skeleton of $K_2$. The map $\alpha$ is a homotopy equivalence $\alpha:\underline H'\fl
\underline K'$ and we may add $2$-cells to $\underline H'$ in order to get a complex $\underline H_2$ and a homotopy equivalence $\alpha:\underline H_2\fl\underline K_2$.

Let $C$ be the cylinder of the map $\alpha:\underline H_2\fl\underline K_2$ and $\underline C_2$ be the complex defined by the cocartesian square:
$$\diagram{\underline K_1\times[0,1]&\hfl{}{}&C\cr\vfl{pr_1}{}&&\vfl{}{}\cr \underline K_1&\hfl{}{}&\underline C_2\cr}$$

The complex $\underline C_2$ contains $\underline K_2$ and $\underline H_2$ as subcomplexes. Moreover these inclusions are homotopy equivalences. Since each loop 
$\widehat\sigma$ is homotopic to the trivial loop in $B$, we may extend the map $\underline K_2\fl B$ to a map $\underline C_2\fl B$ which is trivial on each $S(\sigma)$.
With this new map, we get $B$-spaces $H_2$, $C_2$ and $H'$. Moreover the inclusion $K_1\subset H_2$ belongs to $W(2)$. 

Since the map $\underline H_2\fl B$ is trivial
on $S$, $S(\sigma)$ induces, for each $\sigma\in\Sigma$, a well defined element $\lambda(\sigma)$ in $\pi_2(\underline H'\fl B)=\Gamma_2(H')$.

Since $G=\pi_1(K_1)$ acts on $\Gamma_2(H')$, the correspondence $\lambda$ extends to a morphism of $G$-groups 
$\lambda:F_G[\Sigma]\fl\Gamma_2(H')$.

Let $\Sigma'$ be the set of $2$-cells in $H_2\setminus K_1$. As above, we take an orientation for each cell in $\Sigma'$. Let's take a point in each circle $S(\sigma)$
which is not its base point. So we get a finite set $F\subset S$.

Consider a cell $\sigma\in\Sigma'$. This cell is attached to $H'$ by a map $\psi:(S^2\subset B^2)\fl(\underline H'\fl B)$ corresponding to a commutative diagram:
$$\diagram{S^1&\hfl{}{}&B^2\cr\vfl{\psi_0}{}&&\vfl{\psi_1}{}\cr\underline H'&\hfl{}{}&B\cr}$$
The homotopy class of $\psi$ will be denoted by $u(\sigma)$.

We may modified $\psi$ by homotopy in such a way that $\psi_0$ is transverse to $F$ and $\psi_1$ is trivial in a neighborhood of $\psi_0^{-1}(F)$. By removing to $B^2$ a 
small neighborhood of $\psi_0^{-1}(F)$, we get an element $a$ in $\pi_2(\underline K_1\fl B)=\Gamma_2(K_1)=\Gamma_2(X)$ and the contribution of $\psi$ near 
$\psi_0^{-1}(F)$ is the image under $\lambda$ of some element in $F_G[\Sigma]$. Therefore the class of $\psi$ is on the form $u(\sigma)=a\lambda(u)$ for some 
$a\in\Gamma_2(K_1)=\Gamma_2(X)$ and some $u\in F_G[\Sigma]$.

So there are two maps $v:\Sigma'\fl \Gamma_2(X)$ and $w:\Sigma'\fl F_G[\Sigma]$ such that: $u(\sigma)=v(\sigma)\lambda(w(\sigma))$ for each 
$\sigma\in\Sigma'$.

These two maps induce two morphisms of $G$-groups $F_G[\Sigma']\fl \Gamma_2(X)$ and $F_G[\Sigma']\fl F_G[\Sigma]$ still denoted by $v$ and $w$.

The morphism $w:F_G[\Sigma']\fl F_G[\Sigma]$ induces a morphism $R[\Sigma']\fl R[\Sigma]$ which is the boundary of the cellular complex 
$C_*(\underline H_2,\underline K_1;R)$. But $K_1\subset H_2$ belongs to $W(2)$ and the cellular complex is acyclic. Therefore the morphism $w$ belongs to the class
$W(G,R)$.

Consider an extension $g':H'\fl X$ of $g_1:K_1\fl X$. This extension is determined by the images $\mu(\sigma)$ of $\lambda(\sigma)$ in $\Gamma_2(X)$ and then
by a morphism of $G$-groups $\mu:F_G[\Sigma]\fl \Gamma_2(X)$. Moreover this extension $g'$ has an extension to $H_2$ if and only if we have:
$v(\sigma)\mu(w(\sigma))=1$ for every $\sigma\in\Sigma$.

The correspondence $\sigma\mapsto v(\sigma)^{-1}$ define a morphism of $G$-groups $v':F_G[\Sigma]\fl \Gamma_2(X)$ and $g'$ has an extension to $H_2$ if and
only if we have: $\mu w=v'$. But $\Gamma_2(X)$ is $W(G,R)$-local. Therefore such a $\mu$ exists and $g_1:K_1\fl X$ has an extension $H_2\fl X$.

Since $H_2\subset C_2$ is a homotopy equivalence and $\underline X\fl B$ is a Serre fibration, the extension $H_2\fl X$ has an extension $C_2\fl X$. Hence, by restriction
we get an extension $g_2;K_2\fl X$ of $g_1;K_1\fl X$.
\vskip 12pt
\noi{\bf 2.5 Construction of the maps $g_n:K_n\fl X$ and $g':L'\fl X$.}
\vskip 12pt
Let $M$ be a $W(\Z[G],R)$-local module.

Let $i$ be an integer with $2<i\leq n$. The cellular complex $C_*(\underline K_i,\underline K_{i-1};\Z[G])$ is concentrated in degree $i-1$ and $i$. Moreover this 
complex is $R$-acyclic. Therefore the nonzero differential of this complex belongs to $W(\Z[G],R)$ and the cohomology $\H^*(\underline K_i,\underline K_{i-1};M)$ is
trivial.

Since that is true for $i=3,4,\dots,n$, the cohomology $\H^*(\underline K_n,\underline K_2;M)$ is also trivial. In particular, we have: 
$\H^i(\underline K_n,\underline K_2;\Gamma_i(X))=0$ for every $i$. But these modules are the obstructions in order to have an extension $K_n\fl X$ of $g_2:K_2\fl X$.

Therefore there is no obstruction to obtain such an extension $g_n:K_n\fl X$ and $g_n$ is an extension of $g_1$. But $K_n\fl L'$ is a homotopy equivalence, $K_1\fl L'$ is 
a cofibration and $f$ is a Serre  fibration. Therefore $g_1:K_1\fl X$ has an extension $g':L'\fl X$.

Then we get an extension $L\fl L'\build\fl_{}^{g'}X$ of $K\fl X$ and $X$ is weakly $W(B,R)$-local. Hence $X$ is $W(B,R)$-local and theorem 2 is proven.\cqfd
\vskip 24pt
\noi{\bf 3 Finite localization of modules and $G$-groups.}
\vskip 12pt
Let $A\fl R$ be a homomorphism of rings. Denote by Mod$_A$ the category of right $A$-modules and by $W(A,R)$ the class of morphisms in Mod$_A$ between two finitely
generated free $A$-modules inducing an isomorphism in Mod$_R$.

By formally inverting morphisms in $W(A,R)$ we get the Cohn localization $\Lambda=L(A\fl R)$. The tensor product by $\Lambda$ (where $\Lambda$
is considered as a $A$-bimodule) induces a functor $L^A_R:$ Mod$_A\fl$Mod$_A$ and the morphism $A\fl \Lambda$ induces a morphism $i$ from the identity of Mod$_A$ to 
$L^A_R$.
\vskip 12pt
\noi{\bf 3.1 Theorem:} {\sl The pair $(L^A_R,i)$ is a localization in Mod$_A$ with respect to $W(A,R)$.

  Moreover, for every module $M$ in Mod$_A$ and every finite subset $X\subset L^A_R(M)$, there is a morphism $\alpha:F\fl F'$ in $W(A,R)$ and a commutative diagram in
  Mod$_A$:
  $$\diagram{F&\hfl{\alpha}{}&F'\cr\vfl{}{}&&\vfl{\gamma}{}\cr M&\hfl{i}{}&L^A_R(M)\cr}$$
  such that  $X$ is contained in the image of $\gamma$.}
\vskip 12pt
\noi{\bf Proof:} The fact that $(L^A_R,i)$ is a localization in Mod$_A$ with respect to $W(A,R)$ is well known (see [C1] and [V2] section 1).

Let $X$ be a finite subset in $L^A_R(M)=M\build\otimes_A^{}\Lambda$. This set involves only finitely many elements in $M$ and there is a finitely generated free $A$-module
$P$ and a morphism $f:P\fl M$ such that $X$ is contained in the image of $P\build\otimes_A^{}\Lambda\fl M\build\otimes_A^{}\Lambda$. Therefore there is a morphism
$h:A^n\fl L^A_R(P)$ such that the image of the composite morphism $A^b\fl L^A_R(P)\fl L^A_R(M)$ contains $X$.

Because of Proposition 1.3 in [V2], there is a morphism $\alpha:F\fl F'$ in $W(A,R)$ and morphisms $\lambda:F\fl P$ and $\mu:A^n\fl F'$ such that:
$h=\lambda\alpha^{-1}\mu$. Hence we get a commutative diagram:
$$\diagram{&&A^n\cr &&\vfl{\mu}{}\cr F&\hfl{\alpha}{}&F'\cr\vfl{f\lambda}{}&&\vfl{\gamma}{}\cr M&\hfl{i}{}&L^A_R(M)\cr}$$
with: $\gamma=if\lambda\alpha^{-1}$ and the desired property is proven.\cqfd
\vskip 12pt
Let $G$ be a group and $\Z[G]\fl R$ be a homomorphism of rings. Denote by Gr$_G$ the category of right $G$-groups. We have an inclusion functor
Mod$_{\Z[G]}\subset$ Gr$_G$ and an abelianization functor Gr$_G\fl$ Mod$_G$ inducing a functor Gr$_G\fl$ Mod$_R$. We denote by $W(G,R)$ the class of morphisms in Gr$_G$
between two finitely generated free $G$-groups sent to isomorphisms in Mod$_R$.
\vskip 12pt
\noi{\bf 3.2 Theorem:} {\sl The localization $(L^{\Z[G]}_R,i)$ in Mod$_{\Z[G]}$ with respect to $W(\Z[G],R)$ extends to a localization $(L^G_R,i)$ in Gr$_G$ with
  respect to $W(G,R)$.

  Moreover, for every $G$-group $\Gamma$ and every finite subset $X\subset L^G_R(\Gamma)$, there is a morphism $\alpha:F\fl F'$ in $W(G,R)$ and a commutative diagram in
  Gr$_G$:
  $$\diagram{F&\hfl{\alpha}{}&F'\cr\vfl{\beta}{}&&\vfl{\gamma}{}\cr \Gamma&\hfl{i}{}&L^G_R(\Gamma)\cr}$$
  such that  $X$ is contained in the image of $\gamma$.}
\vskip 12pt
\noi{\bf Proof:} Let $\Omega$ be a set $\{F_0,F_1,F_2,\dots\}$ where each $F_n$ is a free $G$-group of rank $n$. Denote by $W_0(G,R)$ be the set of morphisms in $W(G,R)$ 
with source and target in $\Omega$.

For each $G$-group $\Gamma$, $V(\Gamma)$ is defined to be the
set of diagrams $\Gamma\build\longleftarrow_{}^\beta F\build\fl_{}^\alpha F'$ in Gr$_G$ with $\alpha:F\fl F'$ in $W_0(G,R)$. For each
$v\in V(\Gamma)$, the diagram $v$ will be also denoted by $\Gamma\build\longleftarrow_{}^{\beta_v} F_v\build\fl_{}^{\alpha_v} F'_v$.

The category Gr$_G$ has colimits and the coproduct in Gr$_G$ is the free product $*$. Then, for every subset $J\subset V(\Gamma)$, the free product of the diagrams
$v\in J$ is a diagram: $D(J)=(\Gamma\build\longleftarrow_{}^{\beta_J} F_J\build\fl_{}^{\alpha_J} F'_J)$, and the colimit of this diagram is a $G$-group denoted by
$\Phi(\Gamma,J)$. So we have a cocartesian diagram:
$$\diagram{F_J&\hfl{\alpha_J}{}&F'_J\cr\vfl{\beta_J}{}&&\vfl{\gamma_J}{}\cr\Gamma&\hfl{}{}&\Phi(\Gamma,J)\cr}$$

Denote by $\Phi(\Gamma)$ the $G$-group $\Phi(\Gamma,V(\Gamma))$. The correspondence $\Gamma\mapsto\Phi(\Gamma)$ is clearly a functor $\Phi:$ Gr$_G\fl$ Gr$_G$ and we have a
morphism of functors $j:$ Id$\fl\Phi$.

For each $G$-group $\Gamma$, $V'(\Gamma)$ is defined to be the set of triples $(\alpha:F\fl F',\gamma,\gamma')$ where $\alpha$ is in $W_0(G,R)$
and $\gamma$ and $\gamma'$ are morphisms from $F'$ to $\Gamma$ such that: $\gamma\alpha=\gamma'\alpha$. As above, each $v\in V'(\Gamma)$ will be denoted by
$(\alpha_v:F_v\fl F'_v,\gamma_v,\gamma'_v)$ and, for every subset $J\subset V'(\Gamma)$, the free product of the triples $v\in J$ is a triple
$(\alpha_J:F_J\fl F'_J,\gamma_J,\gamma'_J)$ and the coequalizer of $\gamma_J$, $\gamma'_J:F'_J\fl \Gamma$ is a $G$-group denoted by $\Psi(\Gamma,J)$.

Denote by $\Psi(\Gamma)$ the $G$-group $\Psi(\Gamma,V'(\Gamma))$. As above the correspondence $\Gamma\mapsto\Psi(\Gamma)$ is clearly a functor $\Psi:$ Gr$_G\fl$ Gr$_G$ and
we have a morphism of functors $k:$ Id$\fl\Psi$.
\vskip 12pt
\noi{\bf 3.3 Lemma:} {\sl For every $G$-group $\Gamma$ and every diagram $\Psi\Phi(\Gamma)\build\longleftarrow_{}^\beta F\build\fl_{}^\alpha F'$ in
  $V(\Psi\Phi(\Gamma))$, the morphism $\beta:F\fl \Psi\Phi(\Gamma)$ has an extension $F'\fl \Psi\Phi(\Gamma)$.}
\vskip 12pt
\noi{\bf Proof:} Let us denote $\Phi(\Gamma)$ by $\Gamma_1$ and $\Psi(\Gamma_1)$ by $\Gamma_2$. Since $\Gamma_2$ is a coequalizer, the morphism $k:\Gamma_1\fl\Gamma_2$ is
onto and the morphism $\beta:F\fl \Gamma_2$ has a lifting $\beta_1:F\fl \Gamma_1$.

On the other hand, $\Gamma_1$ is the filtered colimit of the $\Phi(\Gamma,J)$'s with $J$ finite. Then, because $F$ is finitely generated, there is a finite subset
$J\subset V(\Gamma)$ such that $\beta_1:F\fl \Gamma_1$ factors through $\Phi(\Gamma,J)$ and therefore through $F'_J$.

So we have a finite subset $J\subset V(\Gamma)$ and a morphism $\lambda:F\fl F'_J$ such that the composite morphism
$$F\build\fl_{}^\lambda F'_J\fl \Phi(\Gamma,J)\fl \Gamma_1\fl\Gamma_2$$
is the given morphism $\beta:F\fl\Gamma_2$.

Denote by $X$ the $G$-group defined by the cocartesian diagram:
$$\diagram{F&\hfl{\alpha}{}&F'\cr\vfl{\lambda}{}&&\vfl{g}{}\cr F'_J&\hfl{f}{}&X\cr}\leqno{(D)}$$

Let $E$ be a basis of $F$. We have a morphism $\varphi:F_J*F\fl F'_J*F'$ defined by:
$$\forall x\in F_J,\ \ \varphi(x)=\alpha_J(x)$$
$$\forall x\in E,\ \ \varphi(x)=\lambda(x)^{-1}\alpha(x)$$

After abelianization, $\varphi$ is represented by a matrix like:
$$\pmatrix{\alpha_J&-\lambda\cr 0&\alpha}$$
and $\varphi$ is isomorphic to some morphism in $W_0(G,R)$.

Consider the following diagram:
$$\diagram{F_J*F&\hfl{\varphi}{}&F'_J*F'\cr\vfl{\pi}{}&&\vfl{\psi}{}\cr F_J&\hfl{f\alpha_J}{}&X\cr}\leqno{(D')}$$
where $\pi:F_J*F\fl F_J$ is the identity on $F_J$ and the trivial morphism on $F$ and $\psi:F'_J*F'\fl X$ is induced by $f:F'_J\fl X$ and $g:F'\fl X$.

Actually, this diagram is commutative and cocartesian. Therefore $\Gamma\build\longleftarrow_{}^{\beta_J\pi}F_J*F\build\fl_{}^\varphi F'_J*F'$ is isomorphic to a
diagram in $V(\Gamma)$ and there is a morphism $F'_J*F'\fl \Gamma_1$ making the following diagram commute:
$$\diagram{F_J*F&\hfl{\varphi}{}&F'_J*F'\cr\vfl{\beta_J\pi}{}&&\vfl{}{}\cr\Gamma&\hfl{j}{}&\Gamma_1\cr}$$

But the diagram $(D')$ is cocartesian. Then we get a morphism $h:X\fl \Gamma_1$ making the following diagram commute:
$$\diagram{F_J&\hfl{f\alpha_J}{}&X\cr\vfl{\beta_J}{}&&\vfl{h}{}\cr\Gamma&\hfl{j}{}&\Gamma_1\cr}$$

Denote by $\gamma_J$ the composite morphism: $F'_J\fl\Phi(\Gamma,J)\fl \Gamma_1$. So we have:
$$hf\alpha_J=j\beta_J=\gamma_J\alpha_J$$
and $(\alpha_J:F_J\fl F'_J,hf,\gamma_J)$ is isomorphic to some element in $V'(\Gamma)$ and $hf$ and $\gamma_J$ become equal in $\Gamma_2$. So we have a commutative
diagram:
$$\diagram{F'_J&\hfl{f}{}&X\cr\vfl{\gamma_J}{}&&\vfl{kh}{}\cr \Gamma_1&\hfl{k}{}&\Gamma_2\cr}$$
and, because the diagram $(D)$ is commutative, the following diagram is also commutative:
$$\diagram{F&\hfl{\alpha}{}&F'\cr\vfl{\beta_1}{}&&\vfl{hg}{}\cr\Gamma_1&\hfl{k}{}&\Gamma_2\cr}$$
Hence the morphism $\beta=k\beta_1$ has an extension $hg:F'\fl \Gamma_2$.\cqfd
\vskip 12pt
Using functors $\Phi$ and $\Psi$, we have, for every $G$-group $\Gamma$ a sequence:
$$\Gamma\build\fl_{}^j \Gamma_1\build\fl_{}^k \Gamma_2\build\fl_{}^k \Gamma_3\build\fl_{}^k \dots$$
where: $\Gamma_1=\Phi(\Gamma)$, $\Gamma_2=\Psi(\Gamma_1)$, $\Gamma_3=\Psi(\Gamma_2)$, \dots\ . The colimit of this sequence will be denoted by $L^G_R(\Gamma)$.

The correspondence $\Gamma\mapsto L^G_R(\Gamma)$ is clearly a functor $L^G_R$ from Gr$_G$ to itself and the canonical morphism $i:\Gamma\fl L^G_R(\Gamma)$ is a morphism
from the identity of Gr$_G$ to $L^G_R$.
\vskip 12pt
\noi{\bf 3.4 Lemma:} {\sl For every $G$-group $\Gamma$, the $G$-group $L^G_R(\Gamma)$ is $W(G,R)$-local.}
\vskip 12pt
\noi{\bf Proof:} 
Let $\alpha:F\fl F'$ be a morphism in $W(G,R)$ and $\beta:F\fl L^G_R(\Gamma)$ be a morphism of $G$-groups. Since $\Gamma_2\fl L^G_R(\Gamma)$ is onto and $F$ is a free
$G$-group, $\beta:F\fl L^G_R(\Gamma)$ can be lifted through $\Gamma_2$ and, because of lemma 3.3, $\beta$ as an extension $F'\fl L^G_R(\Gamma)$.

Let $\gamma,\gamma':F'\fl L^G_R(\Gamma)$ be two morphisms such that: $\gamma\alpha=\gamma'\alpha=\beta$ and $\gamma_1,\gamma'_1:F'\fl \Gamma_1$ be two liftings of $\gamma$
and $\gamma'$. For each $n>1$, denote by $\gamma_n$ and $\gamma'_n$ the composite morphisms $F'\build\fl_{}^{\gamma_1}\Gamma_1\fl\Gamma_n$ and
$F'\build\fl_{}^{\gamma'_1}\Gamma_1\fl\Gamma_n$. Since $F$ is finitely generated and free, the morphism:
$$\colim\Hom(F,\Gamma_*)\fl\Hom(F,\colim \Gamma_*)$$
is bijective and the two morphisms $\gamma_1\alpha$ and $\gamma'_1\alpha$ are the same in the colimit of $\Hom(F,\Gamma_*)$. Then there is an integer $n>0$ such that
$\gamma_1\alpha$ and $\gamma'_1\alpha$ are the same in $\Hom(F,\Gamma_n)$ and we have: $\gamma_n\alpha=\gamma'_n\alpha$.

Therefore the triple $(\alpha:F\fl F',\gamma_n,\gamma'_n)$ is isomorphic to some triple in $V'(\Gamma_n)$ and we have: $\gamma_{n+1}=\gamma'_{n+1}$.  Hence the two
morphisms $\gamma$ and $\gamma'$ are equal and  $L^G_R(\Gamma)$ is $W(G,R)$-local.\cqfd
\vskip 12pt
\noi{\bf 3.5 Lemma:} {\sl For every $G$-group $\Gamma$, the morphism $i:\Gamma\fl L^G_R(\Gamma)$ is a $W(G,R)$-equivalence.}
\vskip 12pt
\noi{\bf Proof:} Let $\Gamma'$ be a $W(G,R)$-local $G$-group and $f:\Gamma\fl\Gamma'$ be a morphism of $G$-groups.

Let $\Gamma\build\longleftarrow_{}^\beta F\build\fl_{}^\alpha F'$ be a diagram in $V(\Gamma)$. Since $\Gamma'$ is $W(G,R)$-local, the composite
morphism $F\build\fl_{}^\beta\Gamma\build\fl_{}^f\Gamma'$ has a unique extension $F'\fl \Gamma'$. Because this property holds for every diagram in $V(\Gamma)$, the
morphism $F_J\build\fl_{}^{\beta_J}\Gamma\build\fl_{}^f\Gamma'$ has, for every subset $J\subset V(\Gamma)$, a unique extension $F'_J\fl \Gamma'$. Therefore,
for every $J\subset V(\Gamma)$, the morphism $f:\Gamma\fl\Gamma'$ has a unique extension $\Phi(\Gamma,J)\fl \Gamma'$. In particular $f$ has a unique extension $g_1$ from
$\Gamma_1=\Phi(\Gamma)$ to $\Gamma'$.

Suppose we have a unique extension $g_n:\Gamma_n\fl \Gamma'$ of $f$. Let $(\alpha:F\fl F',\gamma,\gamma')$ be a triple in $V'(\Gamma_n)$. We have: $g_n\gamma\alpha=
g_n\gamma'\alpha$ and, because $\Gamma'$ is $W(G,R)$-local, the two morphisms $g_n\gamma$ and $g_n\gamma'$ are equals. Then, because this property is true for every
triple in $V'(\Gamma_n)$, the morphism $g_n$ has a unique extension $g_{n+1}$ from $\Gamma_{n+1}=\Psi(\Gamma_n)$ to $\Gamma'$.

Thus there is, for every $n>0$, a unique extension $g_n:\Gamma_n\fl\Gamma'$ of $f$. Since these extension are compatible, they induce a unique extension 
$g:L^G_R(\Gamma)\fl\Gamma'$ of $f$ and $i:\Gamma\fl L^G_R(\Gamma)$ is a $W(G,R)$-equivalence.\cqfd
\vskip 12pt
Because of these two lemmas, $(L^G_R,i)$ is a localization in Gr$_G$ with respect to the class $W(G,R)$.
\vskip 12pt
\noi{\bf 3.6 Lemma:} {\sl Let $\Gamma$ be a $G$-group and $X$ be a finite subset of $L^G_R(\Gamma)$. Then there is a morphism $\alpha:F\fl F'$ in $W(G,R)$ and a
  commutative diagram of $G$-groups:
  $$\diagram{F&\hfl{\alpha}{}&F'\cr\vfl{\beta}{}&&\vfl{\gamma}{}\cr\Gamma&\hfl{i}{}&L^G_R(\Gamma)\cr}$$
  such that $X$ is contained in the image of $\gamma$.}
\vskip 12pt
\noi{\bf Proof:} Since $\Phi(\Gamma)\fl L^G_R(\Gamma)$ is onto, $X$ can be lifted in a finite subset $X'\in \Phi(\Gamma)$.

For $J_0=V(\Gamma)$, we a cocartesian diagram of $G$-groups:
$$\diagram{F_{J_0}&\hfl{\alpha_{J_0}}{}&F'_{J_0}\cr\vfl{\beta_{J_0}}{}&&\vfl{\gamma_{J_0}}{}\cr \Gamma&\hfl{}{}&\Phi(\Gamma)\cr}$$
But the morphism $\beta_J$ is clearly onto. Therefore $\gamma_{J_0}$ is also onto and $x'$ is contained in the image of $\gamma_{J_0}$. Therefore there is a finite
subset $J$ of $V(\Gamma)$ such that $X'$ is contained in the image of $\gamma_J:F'_J\fl\Phi(\Gamma)$. Thus we have a morphism $\alpha_J:F_J\fl F'_J$ in $W(G,R)$ and a
commutative diagram of $G$-groups:
$$\diagram{F_J&\hfl{\alpha_J}{}&F'_J\cr\vfl{\beta_J}{}&&\vfl{\gamma_J}{}\cr \Gamma&\hfl{}{}&\Phi(\Gamma)\cr}$$
such that $X'$ is contained is the image of $\gamma_J$. The result follows.\cqfd
\vskip 12pt
\noi{\bf 3.7 Lemma:} {\sl Let $\Gamma$ be a commutative $G$-group. Then $L^G_R(\Gamma)$ is also commutative.}
\vskip 12pt
\noi{\bf Proof:} Let $x$ and $y$ be two elements in $L^G_R(\Gamma)$. Because of lemma 3.6, there is a morphism $\alpha:F\fl F'$ in $W(G,R)$ and a commutative diagram of
$G$-groups:
$$\diagram{F&\hfl{\alpha}{}&F'\cr\vfl{\beta}{}&&\vfl{\gamma}{}\cr\Gamma&\hfl{i}{}&L^G_R(\Gamma)\cr}$$
such that $x$ and $y$ are in the image of $\gamma$.

Let $u$ be an element in $F$ and $\gamma':F'\fl L^G_R(\Gamma)$ be the morphism: $v\mapsto \gamma(\alpha(u)v\alpha(u)^{-1})$. For every $u'\in F$ we have:
$$\gamma'\alpha(u')=\gamma(\alpha(u)\alpha(u')\alpha(u)^{-1})=\gamma\alpha(uu'u^{-1})=i\beta(uu'u^{-1})=i\beta(u')=\gamma\alpha(u')$$
because $\Gamma$ is commutative. But $L^G_R(\Gamma)$ is $W(G,R)$-local. Then we have: $\gamma'=\gamma$ and $\gamma$ sends $[\alpha(F),F']$ to $1$.

Let $v$ be an element in $F'$ and $\gamma'':F'\fl L^G_R(\Gamma)$ be the morphism $v'\mapsto \gamma(vv'v^{-1})$. For every $u'\in F$ we have:
$$\gamma''(\alpha(u'))=\gamma(v\alpha(u')v^{-1})=\gamma\alpha(u')$$
because $\gamma$ sends $[\alpha(F),F']$ to $1$. Since $L^G_R(\Gamma)$ is $W(G,R)$-local, we get: $\gamma''=\gamma$ and we have: $\gamma(vv'v^{-1})=\gamma(v')$.

But $x$ and $y$ are in the image of $\gamma$ and they can be lifted in two elements $v$ and $v'$ in $F'$. Therefore we have:
$$\gamma(vv'v^{-1})=\gamma(v')\ \ \Longrightarrow\ \ xyx^{-1}=y\ \ \Longrightarrow\ \ xy=yx$$
and the $G$-group $L^G_R(\Gamma)$ is commutative.\cqfd
\vskip 12pt
\noi{\bf 3.8 Lemma:} {\sl Let $\Gamma$ be a commutative $G$-group. Then $\Gamma\fl L^G_R(\Gamma)$ is a $W(\Z[G],R)$-localization.}
\vskip 12pt
\noi{\bf Proof:} Let $M$ be a commutative $G$-group. Suppose $M$ is $W(G,R)$-local.

Consider a morphism $u:U\fl U'$ in $W(\Z[G],R)$ and a morphism $f:U\fl M$. There exist two free $G$-groups $F$ and $F'$ and a commutative
diagram of $G$-groups:
$$\diagram{F&\hfl{\alpha}{}&F'\cr\vfl{\lambda}{}&&\vfl{\lambda'}{}\cr U&\hfl{u}{}&U'\cr}$$
such that $\lambda$ and $\lambda'$ induce isomorphisms $\H_1(F)\build\fl_{}^\sim U$ and $\H_1(F')\build\fl_{}^\sim U'$. Therefore the morphism $\alpha:F\fl F'$ belongs
to the class $W(G,R)$ and, because $M$ is $W(G,R)$-local, there is a unique morphism $h:F'\fl M$ such that: $f\lambda=h\alpha$. Hence $M$ is $W(\Z[G],R)$-local.

Suppose $M$ is $W(\Z[G],R)$-local. Let $\alpha:F\fl F'$ be a morphism in $W(G,R)$ and $f:F\fl M$ be a morphism of $G$-groups. We have a commutative diagram of $G$-groups:
$$\diagram{F&\hfl{\alpha}{}&F'\cr\vfl{\lambda}{}&&\vfl{\lambda'}{}\cr\H_1(F)&\hfl{u}{}&\H_1(F')\cr}$$
  Since $M$ is commutative, there is a unique morphism $f':\H_1(F)\fl M)$ such that: $f=f'\lambda$.
  
The morphism $u:\H_1(F)\fl\H_1(F')$ belongs clearly to $W(\Z[G],R)$ and $f':\H_1(F)\fl M$ has a unique extension $g':\H_1(F')\fl M$. Therefore the morphism $g=\lambda'g'$
is the unique morphism $F'\fl M$ extending $f:F\fl M$ and $M$ is $W(G,R)$-local.

Hence a commutative $G$-group is $W(G,R)$-local is and only if it is $W(\Z[G],R)$-local.

Let $\Gamma$ be a commutative $G$-group. Because of lemma 3.7, $L^G_R(\Gamma)$ is commutative and then $W(\Z[G],R)$-local.

Let $\Gamma'$ be a $W(\Z[G],R)$-local $G$-module. Then $\Gamma'$ is $W(G,R)$-local and, because $\Gamma\fl L^G_R(\Gamma)$ is a $W(G,R)$-equivalence, the map:
$\Hom(L^G_R(\Gamma),\Gamma')\fl \Hom(\Gamma,\Gamma')$ is bijective. Thus $\Gamma\fl L^G_R(\Gamma)$ is a $W(\Z[G],R)$-localization.

The result follows and theorems 3 and 3.2 are proven.\cqfd
\vskip 12pt
We have few other results about localization of $G$-groups:
\vskip 12pt
\noi{\bf 3.9 Proposition:} {\sl Let $G$ be a group and $\Z[G]\fl R$ be a ring homomorphism. Let $f:\Gamma\fl\Gamma'$ be an epimorphism of $G$-groups. Then the induced
morphism $L^G_R(\Gamma)\fl L^G_R(\Gamma')$ is onto.}
\vskip 12pt
\noi{\bf Proof:} Let $x$ be an element in $L^G_R(\Gamma')$. Because of lemma 3.6, there is a commutative square of $G$-groups:
$$\diagram{F&\hfl{\alpha}{}&F'\cr\vfl{\beta}{}&&\vfl{\gamma}{}\cr \Gamma'&\hfl{}{}&L^G_R(\Gamma')\cr}$$
where $\alpha:F\fl F'$ is in $W(G,R)$ and $x$ is in the image of $\gamma$. Since $F$ is free and $f$ is onto, we can lift $\beta$ through $\Gamma$ and we get a morphism 
$\beta':F\fl \Gamma$ such that: $f\beta'=\beta$. Since $L^G_R(\Gamma)$ is $W(G,R)$-local, the morphism $F\fl \Gamma\fl L^G_R(\Gamma)$ has a unique extension 
$\gamma':F'\fl L^G_R(\Gamma)$ and we have a commutative diagram:
$$\diagram{F&\hfl{\alpha}{}&F'\cr\vfl{\beta'}{}&&\vfl{\gamma'}{}\cr\Gamma&\hfl{}{}&L^G_R(\Gamma)\cr\vfl{f}{}&&\vfl{}{}\cr \Gamma'&\hfl{}{}&L^G_R(\Gamma')\cr}$$
such that the composite morphism $F'\build\fl_{}^{\gamma'} L^G_R(\Gamma)\fl L^G_R(\Gamma')$ is the morphism $\gamma$.
Since $x$ is in the image of $\gamma$, $x$ is also in the image of $L^G_R(\Gamma)\fl L^G_R(\Gamma')$. The result follows.\cqfd
\vskip 12pt
\noi{\bf 3.10 Proposition:} {\sl Let $G$ be a group and $\Z[G]\fl R$ be a ring epimorphism. Let $\Gamma$ be a finitely generated $G$-group. Suppose: 
$\H_1(\Gamma)\build\otimes_G^{}R=0$. Then we have: $L^G_R(\Gamma)=1$.}
\vskip 12pt
\noi{\bf Proof:} Since $\Gamma$ is finitely generated, there is an epimorphism $f:F\fl\Gamma$ where $F$ is a finitely generated free $G$-group. Let $K$ be the kernel of
$f$. We have an exact sequence:
$$\H_1(K)\build\fl_{}^i\H_1(F)\fl \H_1(\Gamma)\fl 0$$
which implies the following exact sequence:
$$\H_1(K)\build\otimes_G^{}R\fl\H_1(F)\build\otimes_G^{}R\fl \H_1(\Gamma)\build\otimes_G^{}R\fl 0$$ 
and $\H_1(K)\build\otimes_G^{}R\build\fl_{}^i\H_1(F)\build\otimes_G^{}R$ is onto. Hence this map has a section and, because $\Z[G]\fl R$ is onto, this section can be
lifted to a morphism $\H_1(F)\fl \H_1(K)$ and then to a morphism $s:F\fl K$. So we have a commutative diagram of $G$-groups:
$$\diagram{F&\hfl{is}{}&F\cr\vfl{s}{}&&\vfl{=}{}\cr K&\hfl{i}{}&F\cr}$$
and the morphism $is:F\fl F$ induces the identity $H_1(F)\build\otimes_G^{}R\fl\H_1(F)\build\otimes_G^{}R$.

Therefore $is$ belongs to $W(G,R)$ and induces an isomorphism $L^G_R(F)\build\fl_{}^\sim L^G_R(F)$. So we have a commutative diagram:
$$\diagram{L^G_R(F)&\hfl{\sim}{}&L^G_R(F)\cr\vfl{}{}&&\vfl{=}{}\cr L^G_R(K)&\hfl{}{}&L^G_R(F)\cr}$$
and $L^G_R(K)\fl L^G_R(F)$ is onto. But the composite map $K\fl F\fl\Gamma$ is trivial and so is $L^G_R(F)\fl L^G_R(\Gamma)$. On the other hand lemma 3.9 implies that 
$L^G_R(F)\fl L^G_R(\Gamma)$ is onto. Therefore $L^G_R(\Gamma)$ is the trivial $G$-group.\cqfd
\vskip 12pt
\noi{\bf 3.11 Proposition:} {\sl Let $G$ be a group, $\Z[G]\fl R$ be a ring homomorphism, $\Gamma$ be a $G$-group and $\Gamma'$ be a $G$-subgroup of $\Gamma$. Then 
$\Gamma$ has a canonical $G\ltimes\Gamma'$-group structure and $\Gamma$ is $W(G,R)$-local if and only if $\Gamma$ is $W(G\ltimes\Gamma',R)$-local.}
\vskip 12pt
\noi{\bf Proof:} The $G\ltimes\Gamma'$ action on $\Gamma$ is given by the following:
$$x.(g,\gamma)=\gamma^{-1}(x.g)\gamma$$
for each $(x,g,\gamma)$ in $\Gamma\times G\times\Gamma'$.

Suppose $\Gamma$ is $W(G\ltimes\Gamma',R)$-local. If $F$ is a finitely generated free $G$-group, there is a finitely generated free $G\ltimes\Gamma'$-group $F^+$ with the
same basis as $F$ and an inclusion $F\subset F^+$ inducing a bijection:
$$\Hom_{G\ltimes\Gamma'}(F^+,\Gamma)\build\fl_{}^\sim\Hom_G(F,\Gamma)$$

Consider a morphism $\alpha:F\fl F'$ in $W(G,R)$. We have a commutative diagram:
$$\diagram{\Hom_{G\ltimes\Gamma'}(F'^+,\Gamma)&\hfl{\alpha^+}{}&\Hom_{G\ltimes\Gamma'}(F^+,\Gamma)\cr\vfl{\sim}{}&&\vfl{\sim}{}\cr
\Hom_G(F',\Gamma)&\hfl{\alpha}{}&\Hom_G(F,\Gamma)\cr}$$

Since $\Gamma$ in $W(G\ltimes\Gamma',R)$-local, $\alpha^+$ is bijective. But vertical morphisms of the diagram are bijective also. Therefore $\alpha$ is a bijection and
$\Gamma$ is $W(G,R)$-local.

Suppose $\Gamma$ in $W(G,R)$-local. Consider a morphism $\alpha:F\fl F'$ in $W(G\ltimes\Gamma',R)$. Let $U$ and $U'$ be two basis of $F$ and $F'$. The morphism $\alpha$ is
defined by the following:
$$\alpha(u)=\build\prod_i^{} f_i(u).(g_i(u),h_i(u))$$
for every $u\in U$, where $f_i$, $g_i$ and $h_i$ are maps from $U$ to $U'$, $G$ and $\Gamma'$.

Let $W$ be the set of all the $h_i(x)$'s. Denote by $\widehat F$ and $\widehat F'$ the $G$-groups freely generated by $U\coprod W$ and $U'\coprod W$. We have a morphism
$\beta:\widehat F\fl\widehat F'$ defined by:
$$\beta(u)=\build\prod_i^{} h_i(u)^{-1}(f_i(u).g_i(u))h_i(u)\hskip 24pt\hbox{and}\hskip 24pt\beta(w)=w$$
for every $u\in U$ and $w\in W$. It is easy to see that $\beta$ is in $W(G,R)$.

Let $f:F\fl \Gamma$ be a morphism of $G\ltimes\Gamma'$-groups. The map $f:U\fl \Gamma$ and the inclusion $\Gamma'\subset\Gamma$ induce a morphism $g:\widehat F\fl \Gamma$.
Since $\Gamma$ is $W(G,R)$-local, there is a unique morphism $g':\widehat F'\fl \Gamma$ such that: $g'\beta=g$. But this equality is equivalent to the following:
$$g(u)=\build\prod_i^{} g'(h_i(u))^{-1}(g'(f_i(u)).g_i(u))g'(h_i(u))\hskip 24pt\hbox{and}\hskip 24pt g'(w)=w$$
for every $u\in U$ and $w\in W$. Then this condition is equivalent to $g'(w)=w$ for each $w\in W$ and:
$$g(u)=\build\prod_i^{} h_i(u)^{-1}(g'(f_i(u)).g_i(u))h_i(x)\leqno{(1)}$$
for every $u\in U$.

Since $\Gamma$ is a $G\ltimes\Gamma'$-group, there is a unique morphism $f':F'\fl \Gamma$ such that: $f'(u')=g'(u')$ for each $u'\in U'$ and the condition $(1)$ is 
equivalent to:
$$f(u)=\build\prod_i^{} h_i(u)^{-1}(f'(f_i(u)).g_i(u))h_i(u)=\build\prod_i^{}f'(f_i(u)).(g_i(u),h_i(u))\ \ \Longleftrightarrow\ \ f(u)=f'\alpha(u))$$ 
and $f':F'\fl\Gamma$ is the unique morphism such that:  $f=f'\alpha$. Hence $\Gamma$ is $W(G\ltimes\Gamma',R)$-local.\cqfd
\vskip 12pt
\noi{\bf 3.12 Proposition:} {\sl Let $G\fl G'$ be a morphism of groups, $\Z[G']\fl R$ be a ring homomorphism and $\Gamma$ be a $G'$-group. Then we have the following:

$\bullet$ if $\Gamma$ is $W(G',R)$-local then $\Gamma$ is $W(G,R)$-local

$\bullet$ if $G\fl G'$ is onto and $\Gamma$ is $W(G,R)$-local then $\Gamma$ is $W(G',R)$-local.}
\vskip 12pt
\noi{\bf Proof:} Suppose $\Gamma$ is $W(G',R)$-local. Consider a morphism $\alpha:F_1\fl F_2$ in $W(G,R)$. Groups $F_1$ and $F_2$ are free $G$-groups generated by two 
basis $U_1$ and $U_2$. Let $F'_1$ and $F'_2$ be the free $G'$-groups generated by $U_1$ and $U_2$. Using the morphism $G\fl G'$, we get a commutative square:
$$\diagram{F_1&\hfl{\alpha}{}&F_2\cr\vfl{}{}&&\vfl{}{}\cr F'_1&\hfl{\alpha'}{}&F'_2\cr}$$
with $\alpha'$ in $W(G',R)$. Thus we get a commutative square:
$$\diagram{\Hom_{G'}(F'_2,\Gamma)&\hfl{\alpha'^*}{}&\Hom_{G'}(F'_1,\Gamma)\cr\vfl{}{}&&\vfl{}{}\cr \Hom_G(F_2,\Gamma)&\hfl{\alpha^*}{}&\Hom_G(F_1,\Gamma)\cr}$$
where vertical arrows are isomorphisms. Since $\Gamma$ is $W(G',R)$-local, the map $\alpha'^*$ is bijective and so is $\alpha^*$. Hence $\Gamma$ is $W(G,R)$-local.

Suppose $G\fl G'$ is onto and $\Gamma$ is $W(G,R)$-local. Let $\alpha':F'_1\fl F'_2$ be a morphism in $W(G',R)$. Groups $F'_1$ and $F'_2$ are free $G'$-groups generated 
by two basis $U_1$ and $U_2$. Let $F_1$ and $F_2$ be the free $G$-groups generated by $U_1$ and $U_2$. Since $G\fl G'$ is onto, the morphism $\alpha'$ can be lifted to
a morphism $\alpha:F_1\fl F_2$ and we have a commutative square:
$$\diagram{F_1&\hfl{\alpha}{}&F_2\cr\vfl{}{}&&\vfl{}{}\cr F'_1&\hfl{\alpha'}{}&F'_2\cr}$$
with $\alpha$ in $W(G,R)$. Thus we get a commutative square:
$$\diagram{\Hom_{G'}(F'_2,\Gamma)&\hfl{\alpha'^*}{}&\Hom_{G'}(F'_1,\Gamma)\cr\vfl{}{}&&\vfl{}{}\cr \Hom_G(F_2,\Gamma)&\hfl{\alpha^*}{}&\Hom_G(F_1,\Gamma)\cr}$$
where vertical arrows are isomorphisms. Since $\Gamma$ is $W(G,R)$-local, the map $\alpha^*$ is bijective and so is $\alpha'^*$. Hence $\Gamma$ is $W(G',R)$-local.\cqfd
\vskip 24pt
\noi{\bf 4 Proof of theorem 4.}
\vskip 12pt
As before $B$ is a path-connected pointed space $B$ with fundamental group $\pi$. From now on, the ring homomorphism $\Z[\pi]\fl R$ will be supposed to be onto.

Let us denote by $\C_B$ the category of pointed connected $B$-spaces $X$ inducing an epimorphism $\pi_1(X)\fl\pi$.

If $G$ is a group and $\alpha:G\fl\pi$ an epimorphism of groups, $\C_B(G)$ will denote the category of spaces $X$ in $\C_B$ together with a group
homomorphism $G\fl\pi_1(X)$ such that the composite morphism $G\fl\pi_1(X)\fl\pi$ is the given morphism $\alpha$.

Notice that for every space $X$ in $\C_B(G)$, $\Gamma_2(X)$ is a $\pi_1(X)$-group and then a $G$-group. So we get a functor from $\C_B(G)$ to the category of $G$-groups.

In the beginning of the proof of theorem 3.2 two functors $\Phi$ and $\Psi$ from Gr$_G$  to itself were defined as well as two morphisms of functors $j:$ Id$\fl\Phi$ and
$k:$ Id$\fl \Psi$ in such a way that the localization functor $L^G_R$ is the colimit of the sequence:
$$\hbox{Id}\build\fl_{}^j\Phi \build\fl_{}^k\Psi\Phi \build\fl_{}^k\Psi^2\Phi \build\fl_{}^k\Psi^3\Phi\build\fl_{}^k\dots$$
\vskip 12pt
\noi{\bf 4.1 Lemma:} {\sl Let $G\fl\pi$ be an epimorphism of groups and $X$ be an object of $\C_B(G)$. Then there exist an object 
$Y$ in $\C_B(G)$, an $W(B,R)$-equivalence $\varphi;X\fl Y$ and a commutative diagram of $G$-groups:
$$\diagram{\Gamma_2(X)&\hfl{\varphi}{}&\Gamma_2(Y)\cr\vfl{=}{}&&\vfl{}{\sim}\cr\Gamma_2(X)&\hfl{j}{}&\Phi(\Gamma_2(X))\cr}$$
where $\Gamma_2(Y)\fl\Phi(\Gamma_2(X))$ is an isomorphism.}
\vskip 12pt
\noi{\bf Proof:} By construction, there are morphisms $\alpha_i:F_i\fl F'_i$ in $W(G,R)$, morphisms $\beta_i:F_i\fl \pi_2(f)$ and a cocartesian diagram of $G$-groups:
$$\diagram{F&\hfl{\alpha}{}&F'\cr\vfl{\beta}{}&&\vfl{\gamma}{}\cr\pi_2(f)&\hfl{}{}&\Phi(\pi_2(f))\cr}$$
where $\alpha:F\fl F'$ is the coproduct of the $\alpha_i$'s and $\beta$ is induced by the $\beta_i$'s. Let us denote by $I$ the set of indices $i$. 

Let $U_i$ be a basis of $F_i$ and $U'_i$ be a basis of $F'_i$. Denote by $U$ and $U'$ the disjoint union of the $U_i$'s and the $U'_i$'s. The sets $U$ and $U'$ are basis 
of $F$ and $F'$.

Let $S$ be a bouquet of pointed circles indexed by $U'$. For each $v\in U'$, the corresponding circle  in $S$ will be denoted by $S_v$. 

Consider the space $\underline X\vee S$. Using the trivial map $S\fl B$, $\underline X\vee S$ becomes a $B$-space denoted by $X\vee S$. So we have an inclusion map
$\varepsilon:X\subset X\vee S$.

For each $v\in U'$, the map $S_v\subset S$ and the trivial homotopy of the composite map $S_u\subset S\fl B$ induces a well defined element $\lambda(v)$ in 
$\pi_2(\underline X\vee S\fl B)=\Gamma_2(X\vee S)$.
The correspondence $v\mapsto\lambda(v)$ induces a map $U'\fl \Gamma_2(X\vee S)$ and a morphism of $G$-groups $\lambda:F'\fl \Gamma_2(X\vee S)$.

For each $u\in U$, we set: $\mu(u)=\varepsilon_*(\beta(u)^{-1})\lambda\alpha(u)\in\Gamma_2(X\vee S)$. So we have a map $\mu:U\fl \Gamma_2(X\vee S)$ and a morphism of 
$G$-groups $\mu:F\fl \Gamma_2(X\vee S)$.

By adding for each $u\in U$ a $2$-cell to $X\vee S$ attached by $\mu(u)$, we get a $B$-space $Y$ and an inclusion $X\vee S\subset Y$.

By construction, $\Gamma_2(Y)$ is generated by $\Gamma_2(X\vee S)$ and the $\lambda(v)$ with the only relations: $\beta(u)=\lambda(\alpha(u))$. So we 
have a cocartesian diagram of $G$-groups:
$$\diagram{F&\hfl{\alpha}{}&F'\cr\vfl{\beta}{}&&\vfl{\gamma}{}\cr\Gamma_2(X)&\hfl{}{}&\Gamma_2(Y)\cr}$$
and $\Gamma_2(Y)$ is isomorphic to the $G$-group $\Phi(\Gamma_2(X))$.

Let $Z$ be a $W(B,R)$-local $B$-space. Since $X\fl Y$ is a cofibration, the map $\Map(Y,Z)\fl\Map(X,Z)$ is a Serre fibration. Let $g:X\fl Z$ be a map and $E(g)$ be the
corresponding fiber of $\Map(Y,Z)\fl\Map(X,Z)$.

Let $i$ be an index in $I$. By doing the same construction as above with the index $i$ (and the basis $U_i$ and $U'_i$) we get a $B$-space $Y_i$. Moreover $Y_i$ is 
obtained  from $X$ by adding only finitely many $1$ and $2$-cells and the cellular chain complex of the pair $(Y_i,X)$ is reduced to the morphism $\alpha_i$. Therefore
the inclusion $X\subset Y_i$ is a $W(B,R)$-equivalence and the fiber $E_i(g)$ of $g$ in the map $\Map(Y_i,Z)\fl \Map(X,Z)$ is weakly contractible. On the other hand
the fiber $E(g)$ is the product of the $E_i(g)$. Hence $E(g)$ is weakly contractible and $X\fl Y$ is a $W(B,R)$-equivalence.\cqfd
\vskip 12pt
\noi{\bf 4.2 Lemma:} {\sl Let $G\fl\pi$ be an epimorphism of groups and $X$ be an object of $\C_B(G)$. Then there exist an object 
$Y$ in $\C_B(G)$, an $W(B,R)$-equivalence $\varphi:X\fl Y$ and a commutative diagram of $G$-groups:
$$\diagram{\Gamma_2(X)&\hfl{\varphi}{}&\Gamma_2(Y)\cr\vfl{=}{}&&\vfl{}{\sim}\cr\Gamma_2(X)&\hfl{k}{}&\Psi(\Gamma_2(X))\cr}$$
where $\Gamma_2(Y)\fl\Psi(\Gamma_2(X))$ is an isomorphism.}
\vskip 12pt
\noi{\bf Proof:} By construction, we have a set $I$ and, for each $i\in I$, we have morphisms $\alpha_i:F_i\fl F'_i$ in $W(G,R)$, morphisms $\gamma_i$ and $\gamma'_i$ 
from $F'_i$ to $\pi_2(f)$ with the following  properties:

$\bullet$ for each $i$, we have: $\gamma_i\alpha_i=\gamma'_i\alpha_i$

$\bullet$ $\Psi(\Gamma_2(X))$ is the coequalizer of the morphisms $\gamma,\gamma':F'\fl\Gamma_2(X)$, where $\alpha:F\fl F'$ is the coproduct of the $\alpha_i$'s and 
$\gamma$ and $\gamma'$ are induced by the $\gamma_i$'s and the $\gamma'_i$'s.

The morphism $\gamma_i\alpha_i=\gamma'_i\alpha_i$ will be denoted by $\beta_i$. So we have a morphism $\beta:F\fl\Gamma_2(X)$ induced by the $\beta_i$'s.

By replacing the map $f:\underline X\fl B$ by a Serre fibration, we get a $B$-space $X'=(\underline X',f')$ and a map $X\fl X'$ of $B$-spaces such that $\underline X\fl
\underline X'$ is a homotopy equivalence and $f'$ a Serre fibration. Let $E$ by the fiber of $f'$. So we have: $\Gamma_2(X)\simeq\Gamma_2(X')\simeq\pi_1(E)$.

Let us take, for each $i\in I$ a basis $U_i$ of $F_i$ and a basis $U'_i$ of $F'_i$.

Let $S_i$ be a bouquet of circles indexed by the set $U_i\times G$ and $S'_i$ be a bouquet of circles indexed by the set $U'_i\times G$. The bouquet, for all $i\in I$,
of the $S_i$'s and the $S'_i$ will be denoted by $S$ and $S'$. The group $G$ acts (on the right) on all these spaces and we have:
$$\pi_1(S_i)\simeq F_i\hskip 24pt\hbox{and}\hskip 24pt \pi_1(S'_i)\simeq F'_i$$
$$\pi_1(S)\simeq F\hskip 24pt\hbox{and}\hskip 24pt \pi_1(S')\simeq F'$$
Moreover morphisms $\alpha_i:F_i\fl F'_i$ and $\alpha:F\fl F'$ may be represented by $G$-maps $\varphi_i:S_i\fl S'_i$ and $\varphi:S\fl S'$.

For every pointed space $Y$ (with base point $*$), we denote by $Y\wedge[0,1]$ the space $Y\times[0,1]/(*\times[0,1])$. So $Y\wedge[0,1]$ is a pointed space.

Let $C_i$ be the cylinder of $\varphi_i:S_i\fl S'_i$ and $C$ be the cylinder of $\varphi:S\fl S'$. The $G$-space $C$ is the bouquet of the $C_i$'s.

We have inclusions $S_i\subset C_i$ and $S\subset C$ and homotopy equivalences $C_i\fl S'_i$ and $C\fl S'$. Let $\Sigma_i$ be the subspace 
$S_i\wedge[0,1]\cup C_i\wedge\{0,1\}$ of $C_i\wedge[0,1]$ and $\Sigma$ be the subspace $S\wedge[0,1]\cup C\wedge\{0,1\}$ of $C\wedge[0,1]$. The $G$-space $\Sigma$ is
the bouquet of the $\Sigma_i$'s.

We have two maps $C_i\fl E$ induced by the morphisms $\gamma_i$ and $\gamma'_i$ from $F'$ to $\pi_1(E)$ and these maps induce a $G$-map $C_i\wedge\{0,1\}\fl E$.
Because of the equality $\gamma_i\alpha_i=\gamma'_i\alpha_i$, this map has an extension $\lambda_i:\Sigma_i\fl E$. Denote by $\lambda:\Sigma\fl E$ the $G$-map induced by
the $\lambda_i$'s.

Let $E'_i$ and $E'$ be the spaces defined by the following cocartesian diagrams:
$$\diagram{\Sigma_i&\hfl{}{}&C_i\wedge[0,1]\cr\vfl{\lambda_i}{}&&\vfl{}{}\cr E&\hfl{}{}&E'_i\cr}\hskip 48pt
\diagram{\Sigma&\hfl{}{}&C\wedge[0,1]\cr\vfl{\lambda}{}&&\vfl{}{}\cr E&\hfl{}{}&E'\cr}$$

Because of Van Kampen theorem, $\pi_1(E'_i)$ is the coequalizer of the maps $\gamma_i$ and $\gamma'_i$ and $\pi_1(E')$ is the coequalizer of $\gamma$ and $\gamma'$.

We define spaces $\underline Y_i$ and $\underline Y$ by push-out:
$$\diagram{E&\hfl{}{}&E'_i\cr\vfl{}{}&&\vfl{}{}\cr \underline X'&\hfl{}{}&\underline Y_i\cr}\hskip 48pt
\diagram{E&\hfl{}{}&E'\cr\vfl{}{}&&\vfl{}{}\cr \underline X'&\hfl{}{}&\underline Y\cr}$$
With the map $f':\underline X'\fl B$ and trivial maps $E'_i\fl B$ and $E'\fl B$, we get $B$-spaces $Y_i$ and $Y$ and maps of $B$-spaces $X'\fl Y_i$ and $X'\fl Y$.

Since $\Gamma_2(Y)$ is the coequalizer of $\gamma$ and $\gamma'$, we have: $\Gamma_2(Y)\simeq \Psi(\Gamma_2(X))$.

On the other hand the chain complex of $(Y_i,X')$ is equivalent to $\H_1(F_i)\build\fl_{}^{\alpha_i}\H_1(F'_i)$ concentrated in degree $2$ and $3$. Since $\alpha_i$
belongs to $W(G,R)$, the map $X'\fl Y_i$ is a $W(B,R)$-equivalence.

In the same way as above, we show that $X'\fl Y$ is also a $W(B,R)$-equivalence. The result follows.\cqfd
\vskip 12pt
\noi{\bf 4.3 Lemma:} {\sl Let $G\fl\pi$ be an epimorphism of groups and $X$ be an object of $\C_B(G)$. Then there exist an object 
$Y$ in $\C_B(G)$ and a $W(B,R)$-equivalence $X\fl Y$ such that the induced morphism of $G$-groups $\Gamma_2(X)\fl\Gamma_2(Y)$ is a $W(G,R)$-localization.}
\vskip 12pt
\noi{\bf Proof:} Let $\Gamma_0$ be the $G$-group $\Gamma_2(X)$. We have a sequence:
$$\Gamma_0\fl\Gamma_1\fl\Gamma_2\fl\Gamma_3\fl\dots$$
where $\Gamma_1=\Phi(\Gamma_0)$ and $\Gamma_{n+1}=\Psi(\Gamma_n)$ for each $n>0$. Because of lemmas 4.1 and 4.2, there are $B$-spaces $X_n$ and $W(B,R)$-equivalences
$$X=X_0\fl X_1\fl X_2\fl X_3\fl\dots$$ 
such that the induced map $\Gamma_2(X_n)\fl\Gamma_2(X_{n+1})$ is, for each $n\geq0$, isomorphic to $\Gamma_n\fl\Gamma_{n+1}$.

Let $Y$ be the colimit of the $X_n$'s. Then $X\fl Y$ is a $W(B,R)$-equivalence and the induced map $\Gamma_2(X)\fl\Gamma_2(Y)$ is a $W(G,R)$-localization.\cqfd
\vskip 12pt
\noi{\bf 4.4 Proof of theorem 4:} Let $X$ be a $B$-space in $\C_B$ and $G=\pi_1(X)$ be the fundamental group of $\underline X$. 

Consider the case $n=2$. The space $X$ is in the category $\C_B(G)$ and, because of lemma 4.3, there is a $B$-space $Y_2$ in $\C_B$ and a 
$W(B,R)$-equivalence $X\fl Y_2$ inducing a $W(G,R)$-localization $\Gamma_2(X)\fl\Gamma_2(Y_2)$. We have a commutative diagram
with exact lines:
$$\diagram{\pi_2(B)&\hfl{}{}&\Gamma_2(X)&\hfl{i}{}&G&\hfl{}{}&\pi&\hfl{}{}&1\cr\vfl{=}{}&&\vfl{\lambda}{}&&\vfl{\mu}{}&&\vfl{=}{}&&\cr
\pi_2(B)&\hfl{}{}&\Gamma_2(Y_2)&\hfl{j}{}&G'&\hfl{}{}&\pi&\hfl{}{}&1\cr}$$
where $G'$ is the group $\pi_1(Y_2)$.

The morphisms $j$ and $\mu$ induce an epimorphism of groups $\varphi:G\ltimes\Gamma_2(Y_2)\fl G'$. Moreover the $G\ltimes\Gamma_2(Y_2)$-action on $\Gamma_2(Y_2)$ induced 
by $\varphi$ corresponds to the given $G$-action and the conjugation action of $\Gamma_2(X)$. Then, because of proposition 3.11, $\Gamma_2(Y_2)$ is 
$W(G\ltimes\Gamma_2(Y_2),R)$-local and, because of proposition 3.12, $\Gamma_2(Y_2)$ is $W(G',R)$-local.

Suppose $n>2$. The $G$-module $\Gamma_n(X)$ is a commutative $G$-group and $\Gamma_n(X)\fl L^G_R(\Gamma_n(X))$ is a $W(\Z[G],R)$-localization. Then we have a sequence of 
$G$-groups:
$$\Gamma_n(X)\fl \Phi(\Gamma_n(X))\fl\Psi\Phi(\Gamma_n(X))\fl\Psi^2\Phi(\Gamma_n(X))\fl\Psi^3\Phi(\Gamma_n(X))\fl\dots$$
and $L^G_R(\Gamma_n(X))$ is the colimit of this sequence. But $L^G_R(\Gamma_n(X))$ is commutative. Therefore we may replace functors $\Phi$ and $\Psi$ by their 
abelianizations $\Phi'$ and $\Psi'$ defined by:
$$\Phi'(\Gamma)=\H_1(\Phi(\Gamma))\hskip 48 pt \Psi'(\Gamma)=\H_1(\Psi(\Gamma))$$
and $L^G_R(\Gamma_n(X))$ is the colimit of the sequence:
$$\Gamma_n(X)\fl \Phi'(\Gamma_n(X))\fl\Psi'\Phi'(\Gamma_n(X))\fl\Psi'^2\Phi'(\Gamma_n(X))\fl\Psi'^3\Phi'(\Gamma_n(X))\fl\dots$$

Then we may apply the same method as above, by replacing $G$-group by $G$-modules and circles by $n-1$-spheres. Then there is a $B$-space $Y_n$ and a 
$W(B,R)$-equivalence $X\fl Y_n$ such that the morphism $\Gamma_i(X)\fl\Gamma_i(Y_n)$ is a bijection if $i<n$ and a $W(\Z[G],R)$-localization if $i=n$. 

By iterating this construction, we may construct $B$-spaces $Y_n$, $Y_{n+1}$, $Y_{n+2}$, \dots and $W(B,R)$-equivalences:
$$Y_n\fl Y_{n+1}\fl Y_{n+2}\fl\dots$$
such that the following holds, for each $p\geq n$:

$\bullet$ $\pi_1(Y_p)\fl\pi_1(Y_{p+1})$ is an isomorphism

$\bullet$ $\Gamma_i(Y_p)\fl\Gamma_i(Y_{p+1})$ is an isomorphism for each $i\leq p$ and a $W(\pi_1(Y_p),R)$-localization for $i=p+1$.

Therefore $\Gamma_i(Y_p)$ is $W(\pi_1(Y_p),R)$-local for every $i\leq p$ and, because of theorem 2, the colimit $Y$ of the $Y_p$'s is $W(B,R)$-local and $X\fl Y$ is 
a $W(B,R)$-localization. The result follows.\cqfd
\vskip 24pt

\noi{\bf References: }

\begin{list}{}{\leftmargin 55pt \labelsep 10pt \labelwidth 40pt \itemsep 0pt}
\item[{[B]}] A. K. Bousfield -- {\sl The localization of spaces with respect to homology}, Topology {\bf 14} (1975), 133--150.
\item[{[CR]}] S. Cappell and D. Ruberman -- {\sl Imbeddings and homology cobordisms of lens spaces}, Comm. Math. Helv. {\bf 63} (1988), 75--89.
\item[{[CS]}] S. Cappell and J. Shaneson -- {\sl The codimension two placement problem and homology equivalent manifolds}, Ann. of Math. {\bf 99} (1974), 277--348.
\item[{[C1]}] P. M. Cohn -- {\sl Free rings and their relations}, Academic Press (1985).
\item[{[C2]}] P. M. Cohn -- {\sl Inversive localization in Noetherian rings}, Comm. Pure. Appl. Math. {\bf 26} (1973), 679–691.
\item[{[D]}] W. G. Dwyer -- {\sl Noncommutative localization in homotopy theory}, Noncommutative localization in algebra and topology, LMS Lecture Note Series {\bf 330},
  Cambridge University Press (2006), 24-39.
\item[{[FV]}] M. Farber and P. Vogel -- {\sl The Cohn localization of the free group ring}, Math. Proc. Camb. Phil. Soc {\bf 111} (1992), 433--443.
\item[{[LD]}] J. Y. Le Dimet -- {\sl Cobordisme d'enlacements de disques}, Bull. Soc. Math. France {\bf 116} (1988).
\item[{[NR]}] A. Neeman and A. Ranicki -- {\sl Noncommutative localization and chain complexes I. Algebraic K- and L-theory}, math.RA/0109118 (2001).
\item[{[Q]}] D. Quillen -- {\sl The spectrum of an equivariant cohomology ring I \& II}, Ann. Math. {\bf 94} (1971), 549--602.
\item[{[R]}] A. Ranicki -- {\sl Noncommutative localization in topology}, math.AT/0303046 (2003).
\item[{[V1]}] P. Vogel -- {\sl On the obstruction group in homology surgery}, Publ. Math. I.H.E.S, {\bf 55}, N$^\circ$ 1 (1982), 165--206.
\item[{[V2]}] P. Vogel -- {\sl Cohn localization of finite group rings}, J. Knot Theory Ramifications, {\bf 26}, N$^\circ2$, 1740017 (2017) [38 pages].

\end{list}
\end{document}